\documentclass{elsart}
\usepackage[metapost]{mfpic}
\usepackage{amsmath,amsfonts,amssymb,color,lineno}
\usepackage{cite}
\DeclareMathOperator{\gr}{gph}
\DeclareMathOperator{\ptrog}{int}
\DeclareMathOperator{\bien}{bd}

\newcommand{\norm}[1]{\left\Vert#1\right\Vert}

\newcommand {\R} {\mathbb R}

\newcommand {\B} {\mathbb B}
\newcommand {\bd} {{\rm bd}\,}
\newcommand {\Int} {{\rm int}\,}
\newcommand {\sd} {\partial}

\def\RHS{right-hand side}

\newtheorem{theorem}{Theorem}
\newtheorem{proposition}{Proposition}
\newtheorem{corollary}{Corollary}
\newtheorem{definition}{Definition}
\newtheorem{remark}{Remark}
\newtheorem{example}{Example}

\begin{document}
\journal{Journal of Mathematical Analysis and Applications}

\begin{frontmatter}
\title{About $[q]$-regularity Properties of Collections of Sets}

\author{Alexander Y. Kruger\corauthref{cor}},
\corauth[cor]{Corresponding author.}
\ead{a.kruger@federation.edu.au}
\ead[url]{http://uob-community.ballarat.edu.au/\~{}akruger}
\author{Nguyen H. Thao}
\ead{hieuthaonguyen@students.federation.edu.au}

\address{Centre for Informatics and Applied Optimisation,\\
School of Science, Information Technology and Engineering,\\
Federation University Australia, POB 663, Ballarat, Vic 3350, Australia}

\begin{abstract}
We examine three primal space local H\"older type regularity properties of finite collections of sets, namely, \emph{$[q]$-semiregularity}, \emph{$[q]$-subregularity}, and \emph{uniform $[q]$-re\-gularity} as well as their quantitative characterizations.
Equivalent metric characterizations of the three mentioned regularity properties as well as a sufficient condition of $[q]$-subregularity in terms of Fr\'echet normals are established.
The relationships between $[q]$-regularity properties of collections of sets and the corresponding regularity properties of set-valued mappings are discussed.
\end{abstract}

\end{frontmatter}
\section{Introduction}

Regularity properties of collections of sets play an important role in variational analysis and optimization, particularly as constraint qualifications in establishing optimality conditions and coderivative/subdifferential calculus and in analyzing convergence of numerical algorithms.

The concept of \emph{linear regularity} was first introduced in \cite{BauBor93,BauBor96} as a key condition in establishing linear convergence rates of sequences generated by the cyclic projection algorithm for finding a point in the intersection of a collection of closed convex sets.
This property has proved to be an important qualification condition in the convergence analysis, optimality conditions, and subdifferential calculus, cf., \cite{BurDeng05,LewPan98,KlaLi99,LiNahSin00,LiNgPon07,ZheWeiYao10,BauBorLi99,BauBorTse00,BakDeuLi05,AusDanThi05}.

Recently, when investigating the extremality, stationarity and regularity properties of collections of sets systematically, several other kinds of regularity were introduced in \cite{Kru05.1} and have been further investigated in \cite{Kru06.1,Kru09.1,KruLop12.1,KruLop12.2, KruTha13.1,KruTha13.2,Pen13}.
The \emph{uniform regularity} is the negation of the \emph{approximate stationarity} property of collections of sets which is the main ingredient in extensions of the \emph{extremal principle} \cite{Mor06.1,Kru03.1,Kru04.1}.
It has also proved to be useful in the convergence analysis \cite{LewLukMal09,Luk12,Luk13,AttBolRedSou10,KruTha13.1}.

The regularity properties of collections of sets are closely related to the well known regularity properties of set-valued mappings such as the \emph{linear openness}, \emph{covering}, \emph{metric regularity}, \emph{Aubin property}, and \emph{calmness}.
The H\"older extensions of these properties also play an important role in variational analysis both in theory and in establishing convergence rates of numerical algorithms, cf.
\cite{BorZhuang88,Kum09,YenYaoKie08,Fra89,FraQui12,AnhKruTha, GayGeoJea11,LiMor12,Hua12.1}.

In this paper which continues \cite{KruTha13.2}, we attempt to extend regularity properties of collections of sets to the H\"older setting and establish their primal and dual space characterizations.
We also discuss their relationships with the corresponding regularity properties of set-valued mappings.

In Section~\ref{URtheory}, we discuss three primal space local H\"older type regularity properties of finite collections of sets, namely, \emph{$[q]$-semiregularity}, \emph{$[q]$-subregularity}, and \emph{uniform $[q]$-regularity} as well as their quantitative characterizations.
The main result of this section -- Theorem~\ref{T2.1} -- gives equivalent metric characterizations of the three mentioned regularity properties.
We also give several examples illustrating these regularity properties.
Section~\ref{S4} is dedicated to dual characterizations of the regularity properties.
In Theorem~\ref{T3.1}~(i), we give a sufficient condition of $[q]$-subregularity in terms of Fr\'echet normals.
In Section~\ref{S5}, we present relationships between $[q]$-regularity properties of collections of sets and the corresponding regularity properties of set-valued mappings.

Our basic notation is standard, cf. \cite{Mor06.1,RocWet98}.
For a normed linear space $X$, its topological dual is denoted $X^*$ while $\langle\cdot,\cdot\rangle$ denotes the bilinear form defining the pairing between the two spaces. The closed unit ball in a normed space is denoted $\mathbb{B}$. $B_\delta(x)$ stands for the closed ball with radius $\delta$ and center $x$.
If not specified otherwise, products of normed spaces will be considered with the maximum
type norms.

The Fr\'echet normal cone to a subset $\Omega\subset X$ at $x\in \Omega$ and the Fr\'echet subdifferential of a function $f:X \to\R_\infty=\R\cup\{+\infty\}$ at a point $x$ with $f(x)<\infty$ are defined, respectively, by
\begin{equation*}
N_{\Omega}(x) = \left\{x^* \in X^* \mid \limsup_{u\to x,\,u\in\Omega\setminus\{x\}} \frac {\langle x^*,u-x \rangle}{\|u-x\|} \leq 0 \right\},
\end{equation*}
\begin{equation*}
    \partial f(x)=\left\{x^*\in X^*\mid \liminf_{\substack{u\to x,\,u\neq x}}\frac{f(u)-f(x)-\langle x^*,u-x\rangle}{\norm{u-x}}\ge 0\right\}.
\end{equation*}

For a given set $\Omega$ in $X$, its interior and boundary are denoted, respectively, $\Int\Omega$ and $\bd\Omega$.
The indicator and distance functions associated with $\Omega$ are defined, respectively, by
$$
\delta_{\Omega}(x)=\left\{
                      \begin{array}{ll}
                        0, & \hbox{if }\; x\in \Omega, \\
                        \infty, & \hbox{if }\; x\in X\setminus \Omega,
                      \end{array}
                    \right.
$$
$$
d(x,\Omega)=\inf_{\omega\in\Omega}\norm{x-\omega},\; \forall x\in X.
$$
\section{$[q]$-regularity properties of collections of sets}\label{URtheory}

In this section, we discuss local $[q]$-regularity properties of finite collections of sets and their primal space characterizations.

In the sequel, $\bold{\Omega}$ stands for a collection $\{\Omega_1,\ldots,\Omega_m\}$ of $m$ $(m\ge 2)$ sets in a normed linear space $X$, $\bar x\in \bigcap_{i=1}^m\Omega_i$, and, if not specified otherwise, $q\in (0,1]$.

\subsection{Definitions}

The next definition introduces several mutually related regularity properties of $\bold{\Omega}$ at $\bar x$.

\begin{definition}\label{HolderUR}
\begin{enumerate}
\item
$\bold{\Omega}$ is $[q]$-semiregular at $\bar x$ if there exist positive numbers $\alpha$ and $\delta$ such that
\begin{equation}\label{HR}
\bigcap_{i=1}^m(\Omega_i-x_i)\  \bigcap\ B_{\rho}(\bar x) \neq \emptyset
\end{equation}
for all $\rho \in (0,\delta)$ and all
$x_i\in{X}$ $(i=1,\ldots,m)$ such that $\max\limits_{1\leq i\leq m}\|x_i\| \le (\alpha\rho)^{\frac{1}{q}}$.
\item
$\bold{\Omega}$ is $[q]$-subregular at $\bar x$ if there exist positive numbers $\alpha$ and $\delta$ such that
\begin{equation}\label{Hlr'}
\bigcap_{i=1}^m\left(\Omega_i+ (\alpha\rho)^{\frac{1}{q}}\mathbb{B}\right)\bigcap B_{\delta}(\bar x)\subseteq
\left(\bigcap_{i=1}^m \Omega_i\right)+\rho\mathbb{B}
\end{equation}
for all $\rho\in (0,\delta)$.
\item
$\bold{\Omega}$ is uniformly $[q]$-regular at $\bar x$ if there exist positive numbers $\alpha$ and $\delta$ such that
\begin{equation}\label{HUR}
\bigcap_{i=1}^m(\Omega_i-\omega_i-x_i) \bigcap(\rho\mathbb{B})\neq \emptyset
\end{equation}
for all $\rho \in (0,\delta)$, $\omega_i \in \Omega_i \cap B_{\delta}(\bar{x})$, and all $x_i\in{X}$ $(i=1,\ldots,m)$ such that $\max\limits_{1\leq i\leq m}\|x_i\| \le (\alpha\rho)^{\frac{1}{q}}$.
\end{enumerate}
\end{definition}
When $q=1$, we will skip ``$[1]$'' in the name of the corresponding property and write simply ``semiregular'', ``subregular'', or ``uniformly regular'', cf. \cite[Definition~3.1]{KruTha13.2}.

\begin{remark}\label{rem1}
Among the three regularity properties in Definition~\ref{HolderUR}, the third one is the strongest.
Indeed, condition \eqref{HR} corresponds to taking $\omega_i=\bar x$ in \eqref{HUR}.
To compare properties (ii) and (iii), it is sufficient to notice that condition \eqref{Hlr'} is equivalent to the following one: for any $x\in B_{\delta}(\bar x)$, $\omega_i \in \Omega_i$,  $x_i\in{X}$ $(i=1,\ldots,m)$ such that $\max\limits_{1\leq i\leq m}\|x_i\| \le (\alpha\rho)^{\frac{1}{q}}$, and $\omega_i+x_i=x$ $(i=1,\ldots,m)$, it holds
\begin{equation*}
\bigcap_{i=1}^m(\Omega_i-x)\bigcap (\rho\mathbb{B})\neq \emptyset.
\end{equation*}
This corresponds to taking $\omega_i+x_i=x$ $(i=1,\ldots,m)$ in \eqref{HUR} (with $x\in X$) and possibly choosing a smaller $\delta>0$.
Hence, (iii) $\Longrightarrow$ (i) and (iii) $\Longrightarrow$ (ii).

Properties (i) and (ii) in Definition~\ref{HolderUR} are in general independent -- see examples in Subsection~\ref{ss2.3}.
\end{remark}

\begin{remark}\label{compare1}
The larger the order $q$ is, the stronger the properties in Definition~\ref{HolderUR} are. \end{remark}
\begin{remark}
When $\bar x\in \ptrog \bigcap_{i=1}^m\Omega_i$, all the properties in Definition~\ref{HolderUR} hold true automatically for any $q\in (0,\infty)$.
\end{remark}

\begin{remark}\label{rem3}
When $\Omega_1=\Omega_2=\ldots=\Omega_m$ and $q\in (0,1]$, property (ii) in Definition~\ref{HolderUR} is trivially satisfied (with $\alpha=\delta=1$).
\end{remark}

Normally, it does not make sense to consider properties (ii) and (iii) in Definition~\ref{HolderUR} when $q>1$.
In the next proposition, we assume temporarily that all properties in Definition~\ref{HolderUR} are defined for all $q>1$.

\begin{proposition}\label{p2.1}
Let the sets $\Omega_i$ $(i=1,\ldots,m)$ be closed and $q>1$.
\begin{enumerate}
\item
$\bold{\Omega}$ is $[q]$-subregular at $\bar x$ \quad$\Leftrightarrow$\quad
$\bold{\Omega}$ is uniformly $[q]$-regular at $\bar x$ \quad$\Leftrightarrow$\quad $\bar x\in \ptrog \bigcap_{i=1}^m\Omega_i$.
\item
If $\bar x\in \ptrog \bigcap_{i=1}^m\Omega_i$, then $\bold{\Omega}$ is $[q]$-semiregular at $\bar x$.
\item
If $\bold{\Omega}$ is $[q]$-semiregular at $\bar x$ and the sets of \emph{primal proximal normals} \cite[Definition~4.28]{Pen13} $N^P_{\Omega_i}(\bar x):=\{u\in X \mid \exists r>0,\; d(\bar x+ru,\Omega_i)=r\|u\|\}$ are nontrivial for all $i=1,\ldots,m$ such that $\bar x\in\bien\Omega_i$, then $\bar x\in \ptrog \bigcap_{i=1}^m\Omega_i$.
\end{enumerate}
\end{proposition}

\textbf{Proof.}
$(i)$ The implications $\bar x\in \ptrog \bigcap_{i=1}^m\Omega_i$ \quad$\Rightarrow$\quad $\bold{\Omega}$ is uniformly $[q]$-regular at $\bar x$ \quad$\Rightarrow$\quad
$\bold{\Omega}$ is $[q]$-subregular at $\bar x$ are obvious.
Next we show that $\bold{\Omega}$ is $[q]$-sub\-regular at $\bar x$ \quad$\Rightarrow$\quad $\bar x\in \ptrog \bigcap_{i=1}^m\Omega_i$.

Suppose $\bar x\notin \ptrog \bigcap_{i=1}^m\Omega_i$ while $\bold{\Omega}$ is $[q]$-sub\-regular at $\bar x$, i.e., there exist numbers $\alpha>0$ and $\delta>0$ such that condition \eqref{Hlr'} holds true
for all $\rho\in (0,\delta)$.
Consider a sequence $x_k\to\bar x$ such that $r_k:=d(x_k,\bigcap_{i=1}^m\Omega_i)>0$ $(k=1,2,\ldots)$.
Then
$$
x_k\in\bigcap_{i=1}^m\Omega_i+r_k(1+r_k)\B \subseteq\bigcap_{i=1}^m \left(\Omega_i+r_k(1+r_k)\B\right)
$$
and $x_k\in B_{\delta}(\bar x)$ for all sufficiently large $k$.
Denote $\rho_k:=\alpha^{-1}(r_k(1+r_k))^q$.
Then $\rho_k<\delta$ for all sufficiently large $k$, and it follows from \eqref{Hlr'} that $x_k\in\bigcap_{i=1}^m\Omega_i+\rho_k\B$.
Hence, $r_k\le\rho_k$, and consequently $\alpha\le r_k^{q-1}(1+r_k)^q$.
Letting $k\to\infty$, we arrive at a contradiction: $0<\alpha\le0$.

(ii) is obvious.

(iii) Suppose $\bar x\notin \ptrog \bigcap_{i=1}^m\Omega_i$ and there exist numbers $\alpha\ge0$ and $\delta>0$ such that condition \eqref{HR} holds true
for all $\rho \in (0,\delta)$ and all
$x_i\in{X}$ $(i=1,\ldots,m)$ such that $\max_{1\leq i\leq m}\|x_i\| \le (\alpha\rho)^{\frac{1}{q}}$.
Then $\bar x\in \bien{\Omega_j}$ for some $j$. Choose a nonzero $u\in N^P_{\Omega_j}(\bar x)$.
Then there exists a number $r>0$ such that $d(\bar x+tu,\Omega_j)=t\|u\|$ for all $t\in[0,r]$ \cite[p. 284]{Pen13}.
Denote $\rho_t:=t\|u\|$ and $x_t:=(\alpha\rho_t)^{\frac{1}{q}}\frac{u}{\|u\|}$.
Then $\rho_t<\delta$ and $(\alpha\rho_t)^{\frac{1}{q}}/\|u\|<r$ for all sufficiently small $t$.
Hence, $d(\bar x,\Omega_j-x_t) =d(\bar x+x_t,\Omega_j) =(\alpha\rho_t)^{\frac{1}{q}}$, and it follows from \eqref{HR} that $(\alpha\rho_t)^{\frac{1}{q}}\le\rho_t$, and consequently $0\le\alpha\le\rho_t^{q-1}$.
Letting $t\downarrow0$, we conclude that $\alpha=0$, i.e., $\bold{\Omega}$ is not $[q]$-semi\-regular at $\bar x$.
\qed

\begin{remark}
Unlike $[q]$-subregularity and $[q]$-uniform regularity, when $\bar x\notin \ptrog \bigcap_{i=1}^m\Omega_i$, the property of $[q]$-semiregularity can be fulfilled with $q>1$ if the assumption of the existence of nontrivial primal proximal normals in Proposition~\ref{p2.1} is not satisfied -- see Example~\ref{e2.4} below.
\end{remark}

The regularity properties in Definition~\ref{HolderUR} can be equivalently defined using the following nonnegative constants which provide quantitative characterizations of these properties:
\begin{gather}\label{01.1}
\theta^q[\bold{\Omega}](\bar{x}):= \liminf_{\rho \downarrow 0} \dfrac{(\theta_{\rho}[\bold{\Omega}](\bar{x}))^q}{\rho},
\\\label{01.2}
\zeta^q[\bold{\Omega}](\bar{x}):= \lim_{\delta\downarrow 0}\inf_{0<\rho<\delta} \dfrac{(\zeta_{\rho,\delta}[\bold{\Omega}](\bar{x}))^q} {\rho},
\\\label{01.3}
\hat{\theta}^q[\bold{\Omega}](\bar{x}):= \liminf_ {\substack{\omega_i\to\bar{x},\,\omega_i\in{\Omega_i}\,(i=1,\ldots,m)\\ \rho\downarrow0}} \dfrac{(\theta_{\rho}[\Omega_1-\omega_1,\ldots,\Omega_m-\omega_m](0))^q}{\rho},
\end{gather}
where, for $\rho>0$ and $\delta>0$,
\begin{gather}\label{02.1}
\theta_{\rho}[\bold{\Omega}](\bar{x}):= \sup\left\{r \ge 0 \mid \bigcap_{i=1}^m (\Omega_i - x_i) \bigcap B_{\rho}(\bar{x}) \neq \emptyset,\; \forall x_i \in r\mathbb{B}\right\},
\\\label{02.2}
\zeta_{\rho,\delta}[\bold{\Omega}](\bar{x}):= \sup\left\{r\ge0 \mid\bigcap_{i=1}^m (\Omega_i+ r\mathbb{B})\bigcap B_{\delta}(\bar x)\subseteq
\bigcap_{i=1}^m \Omega_i+\rho\mathbb{B}\right\}.
\end{gather}
When $q=1$, we will not write superscript $1$ in the denotations \eqref{01.1} -- \eqref{01.3}.

Using the equivalent representation of condition \eqref{Hlr'} in Remark~\ref{rem1}, it is not difficult to check that $\hat{\theta}^q[\bold{\Omega}](\bar{x})\le \min\{\theta^q[\bold{\Omega}](\bar{x}),\zeta^q[\bold{\Omega}](\bar{x})\}$.

The next proposition follows immediately from the definitions.

\begin{proposition}\label{theorem11}
\begin{enumerate}
\item
$\bold{\Omega}$ is $[q]$-semiregular at $\bar x$ if and only if $\theta^q[\bold{\Omega}](\bar{x})>0$. Moreover, $\theta^q[\bold{\Omega}](\bar{x})$ is the exact upper bound of all numbers $\alpha$ such that (\ref{HR}) is satisfied.
\item
$\bold{\Omega}$ is $[q]$-subregular at $\bar x$ if and only if $\zeta^q[\bold{\Omega}](\bar{x})>0$. Moreover, $\zeta^q[\bold{\Omega}](\bar{x})$ is the exact upper bound of all numbers $\alpha$ such that (\ref{Hlr'}) is satisfied.
\item
$\bold{\Omega}$ is uniformly $[q]$-regular at $\bar x$ if and only if $\hat{\theta}^q[\bold{\Omega}](\bar{x})>0$. Moreover, $\hat{\theta}^q[\bold{\Omega}](\bar{x})$ is the exact upper bound of all numbers $\alpha$ such that (\ref{HUR}) is satisfied.
\end{enumerate}
\end{proposition}

\begin{remark}
With $q=1$, properties (i) and (iii) in Definition~\ref{HolderUR} were discussed in \cite{Kru06.1} (see also \cite[Properties (R)$_S$ and (UR)$_S$]{Kru09.1}), while property (ii) was introduced in \cite{KruTha13.2}.
Constants \eqref{01.1}, \eqref{01.3}, and \eqref{02.1} (with $q=1$) can be traced back to \cite{Kru96.2,Kru98.1,Kru00.1,Kru02.1,Kru03.1,Kru04.1,Kru05.1}.
\end{remark}

The equivalent representation of constant \eqref{02.1} given in the next proposition can be useful.

\begin{proposition}\label{pr1}\cite[Proposition 3.8]{KruTha13.2}
For any $\rho>0$,
\begin{gather}\label{02.1+}
\theta_{\rho}[\bold{\Omega}](\bar{x}):= \sup\left\{r \ge 0 \mid r\mathbb{B}^m\subseteq\bigcup_{x\in B_{\rho}(\bar{x})}\prod_{i=1}^m (\Omega_i-x)\right\},
\end{gather}
where $\prod_{i=1}^m(\Omega_i-x) =(\Omega_1-x)\times\ldots\times(\Omega_m-x)$ and $\mathbb{B}^m=\prod_{i=1}^m\mathbb{B}$.
\end{proposition}

From Propositions~\ref{theorem11} and \ref{pr1}, we immediately obtain equivalent representations of $[q]$-se\-mi\-regularity and $[q]$-uniform regularity.
\begin{corollary}
\begin{enumerate}
\item
$\bold{\Omega}$ is $[q]$-semiregular at $\bar x$ if and only if there exist  positive numbers $\alpha$ and $\delta$ such that
\begin{equation}\label{HR2}
(\alpha\rho)^{\frac{1}{q}} \B^m\subseteq \bigcup_{x\in B_{\rho}(\bar{x})}\prod_{i=1}^m (\Omega_i-x)
\end{equation}
for all $\rho\in(0,\delta)$.
Moreover, $\theta^q[\bold{\Omega}](\bar{x})$ is the exact upper bound of all numbers $\alpha$ such that (\ref{HR2}) is satisfied.
\item
$\bold{\Omega}$ is uniformly $[q]$-regular at $\bar x$ if and only if there exist  positive numbers $\alpha$ and $\delta$ such that
\begin{equation}\label{umi2}
(\alpha\rho)^{\frac{1}{q}} \B^m\subseteq \bigcap_{\substack{\omega_i\in\Omega_i\cap B_{\delta}(\bar x)\\ (i=1,\ldots,m)}} \bigcup_{x\in\rho\B}\prod_{i=1}^m (\Omega_i-\omega_i-x)
\end{equation}
for all $\rho\in(0,\delta)$.
Moreover, $\hat\theta^q[\bold{\Omega}](\bar{x})$ is the exact upper bound of all numbers $\alpha$ such that (\ref{umi2}) is satisfied.
\end{enumerate}
\end{corollary}
\subsection{Metric characterizations}

The $[q]$-regularity properties of collections of sets  in Definition~\ref{HolderUR} can also be characterized in metric terms.
The next proposition generalizing \cite[Proposition~3.15]{KruTha13.2} provides equivalent metric representations of constants \eqref{01.1} -- \eqref{01.3}.

\begin{proposition}\label{pr2}
\begin{align}\label{newcon}
\theta^q[\bold{\Omega}](\bar{x})&= \liminf_{\substack{x_i\to 0\; (i=1,\ldots,m)\\ \bar x\notin \bigcap_{i=1}^m(\Omega_i-x_i)}}
\frac{\max_{1\le i\le m}\norm{x_i}^q}{d\Big{(}\bar x,\bigcap_{i=1}^m(\Omega_i-x_i)\Big{)}},
\\\label{errcon'}
\zeta^q[\bold{\Omega}](\bar{x})&= \quad\liminf_{\substack{x\to\bar{x}\\ x\notin \bigcap_{i=1}^m\Omega_i}}
\quad\frac{\max_{1\le i\le m} d^q(x,\Omega_i)}{d\Big{(}x,\bigcap_{i=1}^m\Omega_i\Big{)}}
\\\notag
&=\liminf_{\substack{x\to\bar{x}\\ \omega_i\to\bar x,\,\omega_i\in{\Omega_i}\, (i=1,\ldots,m)\\ x\notin \bigcap_{i=1}^m\Omega_i}}\frac{\max_{1\le i\le m}\norm{\omega_i-x}^q}{d\Big{(}x,\bigcap_{i=1}^m\Omega_i\Big{)}},
\\\label{vartheta}
\hat{\theta}^q[\bold{\Omega}](\bar{x})&= \liminf_{\substack{x\to\bar{x}\\x_i \to 0\;(i=1,\ldots,m)\\x\notin\bigcap_{i=1}^m(\Omega_i-x_i)}} \frac{\max_{1\le i\le m}d^q(x+x_i,\Omega_i)} {d\Big{(}x,\bigcap_{i=1}^m(\Omega_i-x_i)\Big{)}}
\\\notag
&= \liminf_{\substack{x\to\bar{x}\\x_i\to0,\, \omega_i\to\bar x,\,\omega_i\in{\Omega_i}\,(i=1,\ldots,m)\\ x\notin\bigcap_{i=1}^m(\Omega_i-x_i)}} \frac{\max_{1\le i\le m} \norm{x+x_i-\omega_i}^q} {d\Big{(}x,\bigcap_{i=1}^m(\Omega_i-x_i)\Big{)}}.
\end{align}
\end{proposition}

\textbf{Proof.}
\emph{Equality \eqref{newcon}}.
Let $\xi$ stand for the \RHS\ of \eqref{newcon}.
Suppose that $\xi>0$ and fix an arbitrary number $\gamma\in (0,\xi)$.
Then there is a number $\delta>0$ such that
\begin{equation}\label{new1}
\gamma d\left(\bar x,\bigcap_{i=1}^m(\Omega_i-x_i)\right)\le \max_{1\le i\le m}\norm{x_i}^q,\; \forall x_i\in \delta\mathbb{B}\; (i=1,\ldots,m).
\end{equation}
Choose a number $\alpha\in (0,\gamma)$ and set $\delta'=\frac{\delta^q}{\alpha}$.
Then, for any $\rho\in (0,\delta')$ and $x_i\in (\alpha\rho)^{\frac{1}{q}}\mathbb{B}\; (i=1,\ldots,m)$, it holds
$\max_{1\le i\le m}\norm{x_i}\le (\alpha\rho)^{\frac{1}{q}}\le (\alpha\delta')^{\frac{1}{q}}=\delta$.
Hence, \eqref{new1} yields
$$
d\left(\bar x,\bigcap_{i=1}^m(\Omega_i-x_i)\right)\le \frac{1}{\gamma}\max_{1\le i\le m}\norm{x_i}^q\le \frac{\alpha}{\gamma}\rho<\rho.
$$
This implies \eqref{HR} and consequently $\theta^q[\bold{\Omega}](\bar{x})\ge\alpha$.
Taking into account that $\alpha$ can be arbitrarily close to $\xi$, we obtain $\theta^q[\bold{\Omega}](\bar{x})\ge \xi$.

Conversely, suppose that $\theta^q[\bold{\Omega}](\bar{x})>0$ and fix an arbitrary number $\alpha\in (0,\theta^q[\bold{\Omega}](\bar{x}))$.
Then there is a number $\delta>0$ such that \eqref{HR} is satisfied for all $\rho\in (0,\delta)$ and $x_i\in (\alpha\rho)^{\frac{1}{q}}\mathbb{B}\; (i=1,\ldots,m)$.
Choose a positive $\delta'<(\alpha\delta)^{\frac{1}{q}}$.
For any $x_i\in \delta'\mathbb{B}\; (i=1,\ldots,m)$, it holds $\max_{1\le i\le m}\norm{x_i}< (\alpha\delta)^{\frac{1}{q}}$.
Pick up a $\rho\in(0,\delta)$ such that $\max_{1\le i\le m}\norm{x_i}=(\alpha\rho)^{\frac{1}{q}}$.
Then \eqref{HR} yields
$$
\alpha d\left(\bar x,\bigcap_{i=1}^m(\Omega_i-x_i)\right)\le \alpha\rho=\max_{1\le i\le m}\norm{x_i}^q.
$$
This implies $\xi\ge \alpha$. Since $\alpha$ can be arbitrarily close to $\theta^q[\bold{\Omega}](\bar{x})$, we deduce
$\xi\ge \theta^q[\bold{\Omega}](\bar{x})$.

\emph{Equality \eqref{errcon'}}.
Let $\xi$ stand for the \RHS\ of \eqref{errcon'}.
Suppose that $\xi>0$ and fix an arbitrary number $\alpha\in (0,\xi)$. Then there is a number $\delta>0$ such that
\begin{equation*}
\alpha d\left(x,\bigcap_{i=1}^m\Omega_i\right) \le \max_{1\le i\le m} d^q(x,\Omega_i),\; \forall x \in B_{\delta}(\bar{x}).
\end{equation*}
If $x\in\bigcap_{i=1}^m\left(\Omega_i+ (\alpha\rho)^{\frac{1}{q}}\mathbb{B}\right)\bigcap B_{\delta}(\bar x)$ for some $\rho\in (0,\delta)$, then $\max_{1\le i\le m} d^q(x,\Omega_i)\le\alpha\rho$, and consequently $d\left(x,\bigcap_{i=1}^m\Omega_i\right) \le\rho$, i.e., $\zeta_{\rho,\delta}[\bold{\Omega}](\bar{x}) \ge(\alpha\rho)^{\frac{1}{q}}$.
Hence, $\zeta^q[\bold{\Omega}](\bar{x})\ge \alpha$. Since $\alpha$ can be arbitrarily close to $\xi$, we obtain $\zeta^q[\bold{\Omega}](\bar{x})\ge \xi$.

Conversely, suppose that $\zeta^q[\bold{\Omega}](\bar{x})>0$ and fix any $\alpha\in (0,\zeta^q[\bold{\Omega}](\bar{x}))$.
Then there is a number $\delta>0$ such that
\eqref{Hlr'} is satisfied for all $\rho\in (0,\delta)$.
Choose a positive number $\delta'<\min\{(\alpha\delta)^{\frac{1}{q}},\delta\}$.
For any $x \in B_{\delta'}(\bar{x})$, it holds
$$
\max_{1\le i\le m}d(x,\Omega_i)\le \norm{x-\bar x}\le \delta'<(\alpha\delta)^{\frac{1}{q}}.
$$
Choose a $\rho\in(0,\delta)$ such that $\max_{1\le i\le m}d(x,\Omega_i)=(\alpha\rho)^{\frac{1}{q}}$.
Then, by \eqref{Hlr'},
$$
\alpha d\left(x,\bigcap_{i=1}^m\Omega_i\right) \le \alpha\rho =\max_{1\le i\le m}d^q(x,\Omega_i).
$$
Hence, $\alpha\le\xi$.
By letting $\alpha\to\zeta^q[\bold{\Omega}](\bar{x})$, we obtain $\zeta^q[\bold{\Omega}](\bar{x})\le\xi$.

\emph{Equality \eqref{vartheta}}.
Let $\xi$ stand for the \RHS\ of \eqref{vartheta}.
Suppose that $\xi>0$ and fix an arbitrary number $\gamma\in(0,\xi)$.
Then there is a number $\delta>0$ such that
\begin{equation}\label{UMI}
    \gamma d\left(x,\bigcap_{i=1}^m(\Omega_i-x_i)\right) \le \max_{1\le i\le m}d^q(x+x_i,\Omega_i)
\end{equation}
for any $x\in B_{\delta}(\bar x)$ and $x_i \in \delta\mathbb{B}\; (i=1,\ldots,m)$.
Fix any positive number $\alpha<\gamma$ and pick up a positive number $\delta'$ satisfying $\delta'+(\alpha\delta')^{\frac{1}{q}}\le \delta$. Then, for any $\rho \in (0,\delta']$, $\omega_i \in \Omega_i \cap B_{\delta'}(\bar x)$ and $a_i \in (\alpha\rho)^{\frac{1}{q}}\mathbb{B}\; (i=1,\ldots,m)$, it holds
$$
\norm{\omega_i-\bar x +a_i} \le \delta'+(\alpha\rho)^{\frac{1}{q}} \le \delta'+(\alpha\delta')^{\frac{1}{q}}\le \delta.
$$
Applying (\ref{UMI}) with $x=\bar x$ and $x_i=\omega_i-\bar x +a_i$, we get
\begin{align*}
d\left(0,\bigcap_{i=1}^m(\Omega_i-\omega_i-a_i)\right) &\le \gamma^{-1}\max_{1\le i\le m} d^q(\omega_i+a_i,\Omega_i)
\\
&\le \gamma^{-1}\max_{1\le i\le m}\norm{a_i}^q \le \frac{\alpha}{\gamma}\rho<\rho.
\end{align*}
Hence, \eqref{HUR} holds true and consequently $\hat{\theta}^q[\bold{\Omega}](\bar{x}) \ge\alpha$.
Taking into account that $\alpha$ can be arbitrarily close to $\xi$, we obtain
$\hat{\theta}^q[\bold{\Omega}](\bar{x}) \ge \xi$.

Conversely, suppose that $\hat{\theta}^q[\bold{\Omega}](\bar{x}) >0$ and fix an arbitrary number $\alpha\in (0,\hat{\theta}^q[\bold{\Omega}](\bar{x}))$.
Then there is some number $\delta>0$ such that \eqref{HUR} is satisfied for all $\rho \in (0,\delta]$, $\omega_i \in \Omega_i \cap B_{\delta}(\bar x)$ and $a_i \in (\alpha\rho)^{\frac{1}{q}}\mathbb{B}\; (i=1,\ldots,m)$.
We pick up some $\delta'>0$ satisfying
\begin{equation}\label{l2}
(\delta'\alpha+(\delta')^q)^{\frac{1}{q}} +\frac{(\delta')^q}{\alpha}+2\delta'<\delta.
\end{equation}
Now, for $x\in B_{\delta'}(\bar x)$ and $x_i \in \delta'\mathbb{B}\; (i=1,\ldots,m)$, we consider two cases.

\emph{Case 1.} There exists some $j\in \{1,\ldots,m\}$ such that
$$
d(x+x_j,\Omega_j) \ge (\delta'\alpha+(\delta')^q)^{\frac{1}{q}}.
$$

Take $\rho=\frac{(\delta')^q}{\alpha}<\delta$, $\omega_i=\bar x$, $a_i=x_i$ $(i=1,\ldots,m)$.
Then $\norm{a_i}\le \delta'=(\alpha\rho)^{\frac{1}{q}}$.
Applying \eqref{HUR}, we find points
$$
x'' \in \bigcap_{i=1}^m(\Omega_i-\bar x-x_i) \bigcap(\rho\mathbb{B})
$$
and
$$
x':=\bar x +x'' \in \bigcap_{i=1}^m(\Omega_i-x_i) \bigcap B_{\rho}(\bar x).
$$
Hence,
\begin{align*}
d\left(x,\bigcap_{i=1}^m(\Omega_i-x_i)\right) &\le \norm{x-x'} \le \norm{x-\bar x} + \norm{x''}
\\
&\le \delta' + \rho = \frac{1}{\alpha}(\delta'\alpha+(\delta')^q)
\\
&\le \frac{1}{\alpha}\max_{1\le i\le m} d^q(x+x_i,\Omega_i),
\end{align*}
and consequently
\begin{equation}\label{case1}
\alpha d\left(x,\bigcap_{i=1}^m(\Omega_i-x_i)\right) \le \max_{1\le i\le m}d^q(x+x_i,\Omega_i).
\end{equation}

\emph{Case 2.} $\max\limits_{1\le i\le m} d(x+x_i,\Omega_i) < (\delta'\alpha+(\delta')^q)^{\frac{1}{q}}$.

Choose $\omega_i \in \Omega_i\; (i=1,\ldots,m)$ such that
$$
\norm{x+x_i-\omega_i}< (\delta'\alpha+(\delta')^q)^{\frac{1}{q}}.
$$
Then, thanks to \eqref{l2},
$$
\norm{\omega_i-\bar x} \le \norm{\omega_i-x-x_i} +\norm{x_i}+\norm{x-\bar x}
<(\delta'\alpha+(\delta')^q)^{\frac{1}{q}}+2\delta'< \delta.
$$
Setting
$$
a_i:=x+x_i-\omega_i\; (i=1,\ldots,m),\; \rho:= \frac{1}{\alpha}\max_{1\le i \le m}\norm{a_i}^q,
$$
we have
$$
\rho<\frac{\delta'\alpha+(\delta')^q}{\alpha} <\delta,\; \norm{a_i} \le (\alpha\rho)^{\frac{1}{q}}\; (i=1,\ldots,m).
$$
Applying \eqref{HUR} again, we find points
$$
x'' \in \bigcap_{i=1}^m(\Omega_i-x-x_i) \bigcap(\rho\mathbb{B})
$$
and
$$
x':=x +x'' \in \bigcap_{i=1}^m(\Omega_i-x_i) \bigcap B_{\rho}(x).
$$
Hence,
$$
d\left(x,\bigcap_{i=1}^m(\Omega_i-x_i)\right) \le \norm{x-x'} \le \rho
= \frac{1}{\alpha} \max_{1\le i\le m}\norm{x+x_i-\omega_i}^q.
$$
Taking infimum in the \RHS\ of the last inequality over $\omega_i \in \Omega_i$ $(i=1,\ldots,m)$, we again arrive at (\ref{case1}).

From (\ref{case1}) we conclude that $\alpha \le\xi$.
Since $\alpha$ can be arbitrarily close to $\hat{\theta}^q[\bold{\Omega}](\bar{x})$, we deduce
$\hat{\theta}^q[\bold{\Omega}](\bar{x}) \le \xi$.

The second equalities in the representations of $\zeta^q[\bold{\Omega}](\bar{x})$ and $\hat{\theta}^q[\bold{\Omega}](\bar{x})$ are straightforward.
\qed

Propositions~\ref{theorem11} and \ref{pr2} imply equivalent metric characterizations of the $[q]$-regularity properties of collections of sets.

\begin{theorem}\label{T2.1}
\begin{enumerate}
\item
$\bold{\Omega}$ is $[q]$-semiregular at $\bar x$ if and only if it is \emph{metrically $[q]$-semiregular} at $\bar x$, i.e., there exist  positive numbers $\gamma$ and $\delta$ such that
\begin{equation}\label{Hri}
    \gamma d\left(\bar x,\bigcap_{i=1}^m(\Omega_i-x_i)\right) \le \max_{1\le i\le m}\norm{x_i}^q,\; \forall x_i \in \delta\mathbb{B}\; (i=1,\ldots,m).
\end{equation}
Moreover, $\theta^q[\bold{\Omega}](\bar{x})$ is the exact upper bound of all numbers $\gamma$ such that (\ref{Hri}) is satisfied.
\item
$\bold{\Omega}$ is $[q]$-subregular at $\bar x$ if and only if it is \emph{metrically $[q]$-subregular} at $\bar x$, i.e., there exist  positive numbers $\gamma$ and $\delta$ such that
\begin{equation}\label{Hlr}
    \gamma d\left(x,\bigcap_{i=1}^m\Omega_i\right) \le \max_{1\le i\le m} d^q(x,\Omega_i),\; \forall x \in B_{\delta}(\bar{x}).
\end{equation}
Moreover, $\zeta^q[\bold{\Omega}](\bar{x})$ is the exact upper bound of all numbers $\gamma$ such that (\ref{Hlr}) is satisfied.
\item
$\bold{\Omega}$ is uniformly $[q]$-regular at $\bar x$ if and only if it is \emph{metrically uniformly $[q]$-regular} at $\bar x$, i.e., there exist  positive numbers $\gamma$ and $\delta$ such that
\begin{equation}\label{umi}
    \gamma d\left(x,\bigcap_{i=1}^m(\Omega_i-x_i)\right) \le \max_{1\le i\le m}d^q(x+x_i,\Omega_i)
\end{equation}
for any $x \in B_{\delta}(\bar{x})$, $x_i \in \delta\mathbb{B}\; (i=1,\ldots,m)$.
Moreover, $\hat\theta^q[\bold{\Omega}](\bar{x})$ is the exact upper bound of all numbers $\gamma$ such that (\ref{umi}) is satisfied.
\end{enumerate}
\end{theorem}

\begin{remark}
With $q=1$, property (\ref{Hlr}) in the above theorem is known as the local \emph{linear regularity}, \emph{linear coherence}, or \emph{metric inequality} \cite{BauBor93,BauBor96,BurDeng05, LewPan98,KlaLi99,LiNahSin00,LiNgPon07, ZheWeiYao10,BauBorLi99,BauBorTse00, BakDeuLi05,AusDanThi05,Iof89,Iof00_, NgaiThe01,ZheNg08,Pen13}.
It was used as the key condition when establishing linear convergence rates of sequences generated by cyclic projection algorithms and a qualification condition for subdifferential and normal cone calculus formulae.
The stronger property (\ref{umi}) is sometimes referred to as \emph{unform metric inequality} \cite{Kru05.1,Kru06.1,Kru09.1}.
Property (\ref{Hri}) with $q=1$ was investigated in \cite{KruTha13.2}.
\end{remark}
\subsection{Examples}\label{ss2.3}

In this subsection, we give several examples illustrating the discussed above regularity properties.
We consider collections of two sets in $\R^2$ having a common point $\bar x = (0,0)$.
In the figures below (except Figure~\ref{HRF3}), the two sets are coloured cyan and yellow, respectively, while their intersection is coloured green.

Below we give two examples of collections of sets that do not satisfy certain $q$-re\-gularity properties when $q=1$, while the corresponding properties are fulfilled when $q =\frac{1}{2}$.

\begin{figure}[!ht]
\begin{center}
\includegraphics{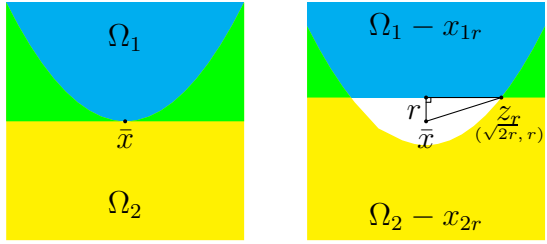}
\end{center}
\caption{Semiregularity vs $[\frac{1}{2}]$-semiregularity}\label{HRF}
\end{figure}

\begin{example}
In the real plane $\R^2$ with the Euclidean norm, consider two sets
$$
\Omega_1:= \left\{(u,v) \in \R^2 \mid v\ge 0\right\},
\quad
\Omega_2:= \left\{(u,v) \in \R^2 \mid v\le u^2\right\},
$$
and the point $\bar x = (0,0) \in \Omega_1 \cap \Omega_2$ (Figure \ref{HRF}).
The collection $\{\Omega_1,\Omega_2\}$ is not semiregular at $\bar x$, while the $\left[\frac{1}{2}\right]$-semiregularity is satisfied at this point.
\end{example}
\textbf{Proof.}
This example is taken from \cite[Figure 8]{Kru09.1}.
We first observe that, for any $r\in(0,1)$ and all $x_1,x_2\in r\mathbb{B}$, it holds
$$
(\Omega_1-x_1) \cap (\Omega_2-x_2) \supseteq (\Omega_1-x_{1r}) \cap (\Omega_2-x_{2r}),
$$
where $x_{1r}=(0,-r)$ and $x_{2r}=(0,r)$.
Besides,
\begin{gather*}
z_{r}:=(\sqrt{2r},r)\in (\Omega_1-x_{1r})\cap(\Omega_2-x_{2r}),
\\
d\left(\bar x,(\Omega_1-x_{1r})\cap (\Omega_2-x_{2r})\right)=\|z_{r}\| =\sqrt{2r+r^2}.
\end{gather*}
Hence, by \eqref{02.1}, for $\rho\in(0,1)$, we have
\begin{gather*}
\theta_{\rho}[\{\Omega_1,\Omega_2\}](\bar{x})= \sup\left\{r \ge 0 \mid \sqrt{2r+r^2}\le\rho\right\} =\sqrt{1+\rho^2}-1,
\end{gather*}
and consequently, by \eqref{01.1},
\begin{gather*}
\theta[\{\Omega_1,\Omega_2\}](\bar{x})= \lim_{\rho\downarrow0} \frac{\sqrt{1+\rho^2}-1}{\rho}=0,
\\
\theta^{\frac{1}{2}}[\{\Omega_1,\Omega_2\}](\bar{x})= \lim_{\rho\downarrow0} \frac{(\sqrt{1+\rho^2}-1)^{\frac{1}{2}}}{\rho} =\frac{1}{\sqrt{2}},
\end{gather*}
which means that
$\{\Omega_1,\Omega_2\}$ is not semiregular at $\bar x$, while it is $\left[\frac{1}{2}\right]$-semiregular at this point.

One can easily show that $\theta_{\rho}[\{\Omega_1-\omega_1,\Omega_2-\omega_2\}](0) \ge \theta_{\rho}[\{\Omega_1,\Omega_2\}](\bar{x})$ for any $\omega_1\in\Omega_1$ and $\omega_2\in\Omega_2$, and consequently, by \eqref{01.3}, $\hat\theta^{\frac{1}{2}}[\{\Omega_1,\Omega_2\}](\bar{x})= \theta^{\frac{1}{2}}[\{\Omega_1,\Omega_2\}](\bar{x})$ and
$\{\Omega_1,\Omega_2\}$ is even $\left[\frac{1}{2}\right]$-uniformly regular at $\bar x$.

Observe also that, for any $x\in\R^2$, $\max_{i=1,2}d(x,\Omega_i)=d(x,\Omega_1\cap\Omega_2)$, and consequently, by \eqref{errcon'}, $\zeta[\{\Omega_1,\Omega_2\}](\bar{x})=1$ and $\{\Omega_1,\Omega_2\}$ is subregular at $\bar x$.
\qed

\begin{figure}[!ht]
\begin{center}
\includegraphics{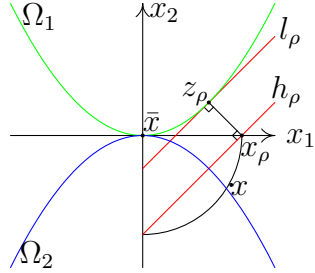}
\end{center}
\caption{Subregularity vs $\left[\frac{1}{2}\right]$-subregularity}\label{HLRF}
\end{figure}

\begin{example}\label{ex2}
In the real plane $\R^2$ with the Euclidean norm, consider two sets
$$
\Omega_1:= \left\{(x,x^2) \in \R^2 \mid x \in \R\right\},
\quad
\Omega_2:= \left\{(x,-x^2) \in \R^2 \mid x \in \R\right\},
$$
and the point $\bar x = (0,0) \in \Omega_1 \cap \Omega_2$ (Figure \ref{HLRF}). The collection $\{\Omega_1,\Omega_2\}$ is not subregular at $\bar x$, while the $\left[\frac{1}{2}\right]$-subregularity is satisfied at this point.
\end{example}
\textbf{Proof.}
We first check that, for each number $\rho\in \left(0,\frac{1}{2}\right)$,
\begin{equation*}
\min\left\{\max_{i=1,2}d(x,\Omega_i) \mid x\in \R^2,\norm{x}=\rho\right\}= d(x_{\rho},\Omega_1)= d(x_{\rho},\Omega_2),
\end{equation*}
where $x_{\rho}:=(\rho,0)$.
By the symmetry of the sets, it suffices to show that
\begin{equation}\label{exa1}
\min\left\{d(x,\Omega_1) \mid x=(x_1,x_2)\in \R^2,\norm{x}=\rho, x_1\ge0, x_2\le0\right\}= d(x_{\rho},\Omega_1).
\end{equation}
Denote $z_{\rho}=(a,a^2):=P_{\Omega_1}(x_{\rho})$ (the metric projection of $x_{\rho}$ onto $\Omega_1$).
Then, with $f(x)=x^2$, we have $f'(z_{\rho})\le 1 =f'\left(\frac{1}{2}\right)$ for any $\rho\in \left(0,\frac{1}{2}\right)$.
Thus, the lines $h_{\rho}$ and $l_{\rho}$ through $x_{\rho}$ and $z_{\rho}$, respectively, with the slope $f'(z_{\rho})$ separate the constraint set in \eqref{exa1} and $\Omega_1$ and consequently, for any $x$ in the constraint set in \eqref{exa1}, it holds
$$
d(x,\Omega_1) \ge d(x,l_{\rho}) \ge d(h_{\rho},l_{\rho})=d(x_{\rho},\Omega_1),
$$
which proves \eqref{exa1}.
One can easily check that $\rho=2a^3+a$ and $d(x_{\rho},z_{\rho})=\sqrt{4a^6+a^4}$.
Hence, by \eqref{errcon'},
\begin{gather*}
\zeta[\{\Omega_1,\Omega_2\}](\bar{x})= \lim_{\rho\downarrow0} \frac{d(x_{\rho},z_{\rho})}{\rho} =\lim_{a\downarrow0} \frac{\sqrt{4a^6+a^4}}{2a^3+a}=0,
\\
\zeta^{\frac{1}{2}}[\{\Omega_1,\Omega_2\}](\bar{x})= \lim_{\rho\downarrow0} \frac{d^{\frac{1}{2}}(x_{\rho},z_{\rho})}{\rho} =\lim_{a\downarrow0} \frac{\sqrt[4]{4a^6+a^4}}{2a^3+a}=1,
\end{gather*}
which means that
$\{\Omega_1,\Omega_2\}$ is not subregular at $\bar x$, while it is $\left[\frac{1}{2}\right]$-subregular at this point.

Observe also that $(\Omega_1-(0,-\varepsilon)) \cap (\Omega_2-(0,\varepsilon))=\emptyset$ for any $\varepsilon>0$.
Hence, by \eqref{02.1} and \eqref{01.1}, $\{\Omega_1,\Omega_2\}$ is not $[q]$-semiregular at $\bar x$ for any $q>0$.
\qed

The above two examples show, in particular, that a collection of sets can be $[q]$-subregular at some point while not being $[q]$-semiregular at this point.
In fact, these two regularity properties are independent.
Next we give an example of a collection of sets that is semiregular at some point while it is not subregular at this point.

\begin{figure}[!ht]
\begin{center}
\includegraphics{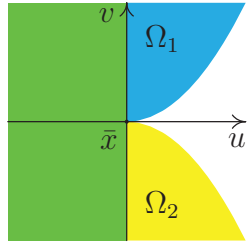}
\end{center}
\caption{Subregularity vs Semiregularity}\label{HRF2}
\end{figure}

\begin{figure}[!ht]
\begin{center}
\includegraphics{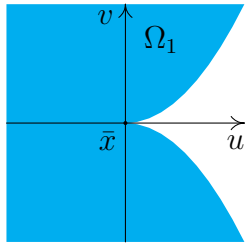}
\end{center}
\caption{$q$-semiregularity}\label{HRF3}
\end{figure}

\begin{example}
In the real plane $\R^2$ with the Euclidean norm, consider two sets
$$
\Omega_1:= \left\{(u,v) \in \R^2 \mid u\le0 \mbox{ or } v\ge u^2\right\},
\;
\Omega_2:= \left\{(u,v) \in \R^2 \mid u\le0 \mbox{ or } v\le -u^2\right\},
$$
and the point $\bar x = (0,0) \in \Omega_1 \cap \Omega_2$ (Figure \ref{HRF2}).
The collection $\{\Omega_1,\Omega_2\}$ is semiregular at $\bar x$, while it is not subregular at this point.
\end{example}
\textbf{Proof.}
The proof of the absence of the subregularity in this example does not differ from that in Example~\ref{ex2}.
Next we show that $\{\Omega_1,\Omega_2\}$ is semiregular at $\bar x$.
For any number $\rho>0$, we set $x_{\rho}:=(-\rho,0)$.
Then
$B_{\rho}(x_{\rho})\subseteq \Omega_i$, i.e., $x_{\rho}+x_i\in\Omega_i$ for any $x_i\in\rho\B$ $(i=1,2)$, and consequently
$$
x_{\rho}\in(\Omega_1-x_1)\cap (\Omega_2-x_2)\cap B_{\rho}(\bar{x}),\; \forall x_i\in\rho\B\;(i=1,2).
$$
Hence, $\theta_{\rho}[\{\Omega_1,\Omega_2\}](\bar{x})\ge\rho$ and $\theta[\{\Omega_1,\Omega_2\}](\bar{x})\ge1$. (One can show that these are actually equalities.)
Thus, $\{\Omega_1,\Omega_2\}$ is semiregular at $\bar x$.
\qed

\begin{example}\label{e2.4}
In the real plane $\R^2$ with the Euclidean norm, consider two sets
$$
\Omega_1:= \left\{(u,v) \in \R^2 \mid u\le0 \mbox{ or } |v|\ge u^2\right\}
$$
(Figure \ref{HRF3}) and $\Omega_2:= \R^2$,
and the point $\bar x = (0,0) \in \Omega_1 \cap \Omega_2$.
The collection $\{\Omega_1,\Omega_2\}$ is $q$-semiregular at $\bar x$ for any $q\in(0,1]$.
\end{example}
\textbf{Proof.}
Obviously $\Omega_2-x=\Omega_2=\R^2$ for any $x\in\R^2$.
Given a $\rho>0$ and an $r\ge0$, using the computations in Example~\ref{ex2}, one can show that $(\Omega_1 - x) \bigcap B_{\rho}(\bar{x}) \neq \emptyset$ for all $x\in r\mathbb{B}$ if and only if $r\le2a^3+a$ where a positive number $a$ satisfies $4a^6+a^4=\rho^2$.
Hence, $\theta_{\rho}[\{\Omega_1,\Omega_2\}](\bar{x})= 2a^3+a$ where $4a^6+a^4=\rho^2$ and consequently
$$
\theta^{q}[\{\Omega_1,\Omega_2\}](\bar{x})= \lim_{a\downarrow0} \frac{(2a^3+a)^q}{a^2\sqrt{4a^2+1}}=+\infty,
$$
i.e., the collection $\{\Omega_1,\Omega_2\}$ is $q$-semiregular at $\bar x$ for any $q\in(0,1]$.

Note that in fact the $q$-semiregularity condition is satisfied for any $q\le2$.
\qed
\section{Dual characterizations}~\label{S4}

This section discusses dual characterizations of $[q]$-regularity properties ($q\in (0,1]$) of a collection of sets $\bold{\Omega}:=\{\Omega_1,\ldots,\Omega_m\}$ at $\bar x\in \bigcap_{i=1}^m\Omega_i$.
We are going to use the notation $\widehat\Omega:= \Omega_1\times\ldots\times\Omega_m\subset X^m$.

Recall that the (normalized) \emph{duality mapping} \cite[Definition~3.2.6]{Lucc06} $J$ between a normed space $Y$ and its dual $Y^*$ is defined as
$$
J(y):=\left\{y^*\in \mathbb{S}_{Y^*}\mid \langle y^*,y\rangle=\norm{y}\right\},\; \forall y\in Y.
$$
Note that $J(-y)=-J(y)$.

The following simple fact of convex analysis is well known (cf., e.g., \cite[Corollary~2.4.16]{Zal02}).

\begin{lem}\label{ll01}
Let $(Y,\|\cdot\|)$ be a normed space.
\begin{enumerate}
\item
$\sd\|\cdot\|(y)=J(y)$ for any $y\ne0$.
\item
$\sd\|\cdot\|(0)=\B^*$.
\end{enumerate}
\end{lem}

Making use of the convention that the topology in $X^m$ is defined by the maximum type norm,
it is not difficult to establish a representation of the duality mapping on $X^m$ (cf. \cite[Proposition~4.2]{KruTha13.2}).
\begin{proposition}\label{prJ}
For each $(x_1,\ldots,x_m)\in X^m$,
\begin{multline*}
J(x_1,\ldots,x_m)= \Biggl\{(x_1^*,\ldots,x_m^*)\in(X^*)^m\mid \sum_{i=1}^m\|x_i^*\|=1;\mbox{ either } x_i^*=0\\
\mbox{or }\left(\|x_i\|=\max_{1\le j\le m} \|x_j\|,\; x_i^*\in\|x_i^*\| J(x_i)\right)\; (i=1,\ldots,m)\Biggr\}.
\end{multline*}
\end{proposition}

In this section, along with the maximum type norm on $X^{m+1}=X\times X^m$, we are going to use another one depending on a parameter $\rho>0$ and defined as follows:
\begin{equation}\label{norm}
\norm{(x,\hat x)}_{\rho} :=\max\left\{\norm{x},\rho \norm{\hat x}\right\},\quad x\in X,\;\hat x\in X^m.
\end{equation}
It is easy to check that the corresponding dual norm has the following representation:
\begin{equation}\label{normd}
\norm{(x^*,\hat x^*)}_{\rho} =\|x^*\|+\rho^{-1}\|\hat x^*\|,\quad x^*\in X^*,\;\hat x^*\in (X^m)^*.
\end{equation}
Note that if, in \eqref{norm} and \eqref{normd}, $\hat x= (x_1,\ldots,x_m)$ and $\hat x^*= (x_1^*,\ldots,x_m^*)$ with $x_i\in X$ and $x_i^*\in X^*$ ($i=1,2,\ldots,m$), then $\|\hat x\|=\max_{1\le i\le m}\|x_i\|$ and $\|\hat x^*\|=\sum_{i=1}^m\|x_i^*\|$.

The next few facts of subdifferential calculus are used in the proof of the main theorem below.

\begin{lem}[\cite{KruTha13.2}, Lemma~4.3]\label{ll02}
Let $X$ be a normed space and
$\varphi(u,\hat u)=\|(u-u_1,\ldots,u-u_m)\|$ ($u\in X$, $\hat u:=(u_1,\ldots,u_m)\in X^m$).
Suppose $x\in X$, $\hat x:=(x_1,\ldots,x_m)\in X^m$, and $\hat v:=(x-x_1,\ldots,x-x_m)\ne0$. Then
\begin{align*}
\sd\varphi(x,\hat x)\subseteq\bigl\{\left(x^*,\hat x^*\right)\in X^*\times (X^*)^m \mid &
-\hat x^*\in J(\hat v),
\\&
\hat x^* =(x_1^*,\ldots,x_m^*),\;
x^*=-(x_{1}^*+\ldots+x_{m}^*)\bigr\}.
\end{align*}
\end{lem}

\begin{lem}\label{ll021}
Let $X$ be a normed space, $\varphi:X\to\R_\infty$, $q>0$, and $f(u):=(\varphi(u))^q$ $(u\in X)$.
If $x\in X$ and $\varphi(x)\ne0$, then $\sd f(x)=q(\varphi(x))^{q-1}\sd\varphi(x)$.
\end{lem}
\textbf{Proof}
follows from the standard chain rule for Fr\'echet subdifferentials, cf., e. g., \cite[Corollary 1.14.1]{Kru03.1}.
\qed

\begin{lem}\label{ll03}
Let $X$ be a normed space and $\hat\omega:=(\omega_1,\ldots,\omega_m)\in\widehat\Omega$.
Then $N_{\widehat\Omega}(\hat\omega) =N_{\Omega_1}(\omega_1)\times\ldots\times N_{\Omega_m}(\omega_m)$.
\end{lem}
\textbf{Proof}
follows directly from the definition of the Fr\'echet normal cone.
\qed


The proof of the main theorem of this section relies heavily on two fundamental results of variational analysis: the \emph{Ekeland variational principle} (Ekeland \cite{Eke74}; cf., e.g., \cite[Theorem~2.1]{Kru03.1}, \cite[Theorem 2.26]{Mor06.1}) and  the \emph{fuzzy} (\emph{approximate}) \emph{sum rule} (Fabian \cite{Fab89}; cf., e.g., \cite[Rule~2.2]{Kru03.1}, \cite[Theorem 2.33]{Mor06.1}).
Below we provide these results for completeness.

\begin{lem}[Ekeland variational principle] \label{l01}
Suppose $X$ is a complete metric space, and $f: X\to\R_\infty$ is lower semicontinuous and bounded from below, $\varepsilon>0, \lambda>0$. If
$$
f(v)<\inf_X f + \varepsilon,
$$
then there exists $x\in X$ such that

(a) $d(x,v)<\lambda $,

(b) $f(x)\le f(v)$,

(c) $f(u)+(\varepsilon/\lambda)d(u,x)\ge f(x)$ for all $u\in X$.
\end{lem}

\begin{lem}[Fuzzy sum rule] \label{l02}
Suppose $X$ is Asplund,
$f_1: X \to\R$ is Lipschitz continuous and
$f_2: X \to\R_\infty$
is lower semicontinuous in a neighborhood of $\bar x$ with $f_2(\bar x)<\infty$.
Then, for any $\varepsilon>0$, there exist $x_1,x_2\in X$ with $\|x_i-\bar x\|<\varepsilon$, $|f_i(x_i)-f_i(\bar x)|<\varepsilon$ $(i=1,2)$ such that
$$
\partial (f_1+f_2) (\bar x) \subset \partial f_1(x_1) +\partial f_2(x_2) + \varepsilon\B^\ast.
$$
\end{lem}

The next theorem gives dual sufficient conditions for $[q]$-regularity of collections of sets in Asplund spaces.
Recall that a Banach space is called \emph{Asplund} if any continuous convex function defined on a nonempty open convex set is Fr\'echet differentiable on a dense subset of its domain.
Asplund spaces form a broad subclass of Banach spaces including, e. g., all spaces which admit Fr\'echet differentiable re-norms (in particular, Fr\'echet smooth spaces).
Reflexive spaces are examples of Fr\'echet smooth spaces.
Asplund property of a Banach space is necessary and sufficient for the fulfillment of some basic results involving Fr\'echet normals and subdifferentials (cf. \cite{Kru03.1,Mor06.1}).
See \cite {Phelps} for various properties and characterizations of Asplund spaces.

\begin{theorem}\label{T3.1}
Let $X$ be an Asplund space and $\Omega_1$, \ldots $\Omega_m$ be closed.
\begin{enumerate}
\item
$\bold{\Omega}$ is $[q]$-subregular at $\bar x$ if there exist positive numbers $\alpha$ and $\delta$ such that, for any $\rho\in (0,\delta)$, $x\in B_{\rho}(\bar x)$, $\omega_i\in\Omega_i\cap B_{\rho}(x)$ $(i=1,\ldots,m)$ with $\omega_j\neq x$ for some $j\in\{1,\ldots,m\}$, there is an $\varepsilon>0$ such that, for 
any $x'\in B_{\varepsilon}(x)$, $\hat\omega'_i\in\Omega_i\cap B_{\varepsilon}(\omega_i)$, $x_{i}^*\in N_{\Omega_i}(\omega_i') +\rho\B^*$ $(i=1,\ldots,m)$
satisfying
$\hat v:=(\omega_1'-x',\ldots,\omega_m'-x')\ne0$ and
\begin{gather*}
x_i^*=0
\quad\mbox{if}\quad
\norm{x'-\omega_i'}<\max_{1\le j\le m} \norm{x'-\omega_j'},
\\
\langle x_i^*,x'-\omega_i'\rangle
\ge\|x_i^*\|(\|x'-\omega_i'\|-\varepsilon),
\\
\sum_{i=1}^m\|x_i^*\|={q}\norm{\hat v}^{q-1},
\end{gather*}
it holds
\begin{gather}\label{>al}
\left\|\sum\limits_{i=1}^m{x}_i^*\right\|>\alpha.
\end{gather}
\item
$\bold{\Omega}$ is uniformly $[q]$-regular at $\bar x$ if there are positive numbers $\alpha$ and $\delta$ such that \eqref{>al} holds true
for all $\omega_i\in\Omega_i\cap{B}_\delta(\bar{x})$ and
$x_i^*\in{N}_{\Omega_i}(\omega_i)$ $(i=1,\ldots,m)$ satisfying $\sum_{i=1}^m\|x_i^*\|=1$.
The inverse implication holds true when $q=1$.
\end{enumerate}
\end{theorem}

The proof of Theorem~\ref{T3.1}~(i) consists of a series of propositions providing lower estimates for constant \eqref{errcon'} and, thus, sufficient conditions for $[q]$-subregularity of $\bold{\Omega}$ which can be of independent interest.
Observe that constant \eqref{errcon'} can be rewritten as
\begin{equation}\label{errcon}
\zeta^q[\bold{\Omega}](\bar{x})= \liminf_{\substack{x\to\bar{x},\, \omega_i\to\bar x\; (i=1,\ldots,m)\\ \hat\omega:=(\omega_1,\ldots,\omega_m)\\x\notin \bigcap_{i=1}^m\Omega_i}} \frac{f_q(x,\hat\omega)} {d\left(x,\bigcap_{i=1}^m\Omega_i\right)}
\end{equation}
with function $f_q:X^{m+1}\to{\R}_\infty :={\R}\cup\{+\infty\}$ defined as
\begin{equation}\label{f}
f_q(x,\hat x)=\max_{1\le i\le m} \norm{x-x_i}^q + \delta_{\widehat\Omega}(\hat x),\quad x\in X,\; \hat x:=(x_1,\ldots,x_m)\in X^m,
\end{equation}
where $\delta_{\widehat\Omega}$ is the indicator function of $\widehat\Omega$: $\delta_{\widehat\Omega}(\hat x)=0$ if $\hat x\in\widehat\Omega$ and $\delta_{\widehat\Omega}(\hat x)=+\infty$ otherwise.

\begin{prop}\label{T3.1.1}
Let $X$ be a Banach space and $\Omega_1$, \ldots, $\Omega_m$ be closed.
\begin{enumerate}
\item
$\hat\zeta^q[\bold{\Omega}](\bar x)\le \zeta^q[\bold{\Omega}](\bar{x})$, where
\begin{equation}\label{strslo}
\hat\zeta^q[\bold{\Omega}](\bar x):= \lim_{\rho\downarrow0} \inf_{\substack{\norm{x-\bar x}<\rho\\ \hat\omega=(\omega_1,\ldots,\omega_m)\in \widehat\Omega\\ 0<\max\limits_{1\le i\le m} \norm{x-\omega_i}<\rho}}
\zeta_{\rho}^q[\bold{\Omega}](x,\hat\omega)
\end{equation}
and, for $x\in X$ and $\hat\omega=(\omega_1,\ldots,\omega_m)\in \widehat\Omega$,
\begin{equation}\label{rhoslo}
\zeta_{\rho}^q[\bold{\Omega}](x,\hat\omega):= \limsup_{\substack{(u,\hat v)\to(x,\hat\omega)\\(u,\hat v)\neq (x,\hat\omega)\\ \hat v=(v_1,\ldots,v_m)\in \widehat\Omega}} \frac{\left(\max\limits_{1\le i\le m} \norm{x-\omega_i}^q-\max\limits_{1\le i\le m} \norm{u-v_i}^q\right)_+}{\norm{(u,\hat v)-(x,\hat\omega)}_{\rho}}.
\end{equation}
\item
If $\hat\zeta^q[\bold{\Omega}](\bar x)>0$, then $\bold{\Omega}$ is $[q]$-subregular at $\bar x$.
\end{enumerate}
\end{prop}

\textbf{Proof.}
(i)
Let $\zeta^q[\bold{\Omega}](\bar{x})<\alpha<\infty$. Choose a $\rho\in(0,1)$ and set
\begin{equation}\label{28}
\eta:= \min\left\{\frac{\rho}{2},\frac{\rho}{\alpha},\rho^{\frac{2}{\rho}}\right\}.
\end{equation}
By \eqref{errcon}, there are $x'\in B_{\eta}(\bar x)$ and $\hat\omega'=(\omega_1',\ldots,\omega_m')\in \widehat\Omega$ such that
\begin{equation}\label{29}
   0<f_q(x',\hat\omega')<\alpha d\left(x',\bigcap_{i=1}^m\Omega_i\right).
\end{equation}
Denote $\varepsilon:=f_q(x',\hat\omega')$ and $\mu:=d\left(x',\bigcap_{i=1}^m\Omega_i\right)$. Then $\mu\le \norm{x'-\bar x}\le \eta\le \frac{\rho}{2}<1$.
Observe that $f_q$ is lower semicontinuous.
Applying to $f_q$ Lemma~\ref{l01} with $\varepsilon$ as above and
\begin{equation}\label{30}
\lambda:=\mu(1-\mu^{\frac{\rho}{2-\rho}}),
\end{equation}
we find points $x\in X$ and $\hat\omega=(\omega_1,\ldots,\omega_m)\in X^{m}$ such that
\begin{equation}\label{31}
   \norm{(x,\hat\omega)-(x',\hat\omega')}_{\rho}<\lambda,\; f_q(x,\hat\omega)\le f_q(x',\hat\omega'),
\end{equation}
and
\begin{equation}\label{32}
f_q(u,\hat v)+\frac{\varepsilon}{\lambda}\norm{(u,\hat v)-(x,\hat\omega)}_{\rho} \ge f_q(x,\hat\omega),
\end{equation}
for all $(u,\hat v)\in X\times X^{m}$.
Thanks to \eqref{31}, \eqref{30}, \eqref{28}, and \eqref{29}, we have
$$
\norm{x-x'}<\lambda <\mu \le \norm{x'-\bar x},
$$
\begin{equation}\label{33}
   d\left(x,\bigcap_{i=1}^m\Omega_i\right)\ge d\left(x',\bigcap_{i=1}^m\Omega_i\right)- \norm{x-x'}\ge \mu-\lambda =\mu^{\frac{2}{2-\rho}},
\end{equation}
\begin{equation}\label{34}
\norm{x-\bar x}\le \norm{x-x'}+\norm{x'-\bar x}<2\norm{x'-\bar x}\le 2\eta\le\rho,
\end{equation}
\begin{equation}\label{35}
f_q(x,\hat\omega)\le f_q(x',\hat\omega')<\alpha\mu\le \alpha\eta\le\rho.
\end{equation}
It follows from \eqref{33}, \eqref{34}, and \eqref{35} that
$$
\norm{x-\bar x}<\rho,\; \hat\omega\in \widehat\Omega,\; 0<\max\limits_{1\le i\le m} \norm{x-\omega_i}^q<\rho.
$$
Observe that $\mu^{\frac{\rho}{2-\rho}}\le \eta^{\frac{\rho}{2-\rho}}< \eta^{\frac{\rho}{2}}\le \rho$, and consequently, by \eqref{29} and \eqref{30},
$$
\frac{\varepsilon}{\lambda} <\frac{\alpha\mu}{\lambda} =\frac{\alpha}{1-\mu^{\frac{\rho}{2-\rho}}} <\frac{\alpha}{1-\rho}.
$$
Thanks to \eqref{32} and \eqref{f}, we have
$$
\max\limits_{1\le i\le m} \norm{x-\omega_i}^q-\max\limits_{1\le i\le m} \norm{u-v_i}^q\le \frac{\alpha}{1-\rho}\norm{(u,\hat v)-(x,\hat\omega)}_{\rho}
$$
for all $u\in X$ and $\hat v=(v_1,\ldots,v_m)\in \widehat\Omega$.
It follows that $\zeta_{\rho}^q[\bold{\Omega}](x,\hat\omega)\le \dfrac{\alpha}{1-\rho}$ and consequently
$$
\inf_{\substack{\norm{x-\bar x}<\rho\\\hat\omega=(\omega_1,\ldots,\omega_m)\in \widehat\Omega\\ 0<\max\limits_{1\le i\le m} \norm{x-\omega_i}<\rho}}
\zeta_{\rho}^q[\bold{\Omega}](x,\hat\omega)\le \frac{\alpha}{1-\rho}.
$$
Taking limits in the last inequality as $\rho\downarrow 0$ and $\alpha\to \zeta^q[\bold{\Omega}](\bar{x})$ yields the claimed inequality.

(ii) follows from (i) and Proposition~\ref{theorem11}~(ii).
\qed

\begin{prop}\label{T3.1.2}
Let $X$ be an Asplund space and $\Omega_1$, \ldots, $\Omega_m$ be closed.
\begin{enumerate}
\item
$\hat\zeta^{q*}_1[\bold{\Omega}](\bar x)\le \hat\zeta^q[\bold{\Omega}](\bar x)$, where $\hat\zeta^q[\bold{\Omega}](\bar x)$ is given by \eqref{strslo},
\begin{equation}\label{strsubslo'}
\hat\zeta^{q*}_1[\bold{\Omega}](\bar x):= \lim_{\rho\downarrow 0} \inf_{\substack{\norm{x-\bar x}<\rho\\ \hat\omega=(\omega_1,\ldots,\omega_m)\in \widehat\Omega\\ 0<\max\limits_{1\le i\le m} \norm{x-\omega_i}<\rho}}
\zeta_{\rho,1}^{q*}[\bold{\Omega}](x,\hat\omega)
\end{equation}
and, for $x\in X$ and $\hat\omega=(\omega_1,\ldots,\omega_m)\in \widehat\Omega$,
\begin{equation}\label{subrhoslo}
\zeta_{\rho,1}^{q*}[\bold{\Omega}](x,\hat\omega):= \inf_{\substack{(x^*,\hat y^*)\in\partial f_q(x,\hat\omega)\\ \norm{\hat y^*}<\rho}}\norm{x^*}
\end{equation}
(with the convention that the infimum over the empty set equals $+\infty$).
\item
If $\hat\zeta^{q*}_1[\bold{\Omega}](\bar x)>0$, then $\bold{\Omega}$ is $[q]$-subregular at $\bar x$.
\end{enumerate}
\end{prop}

\textbf{Proof.}
(i) Let $\hat\zeta^q[\bold{\Omega}](\bar x)< \alpha<\infty$.
Choose a $\beta\in (\hat\zeta^q[\bold{\Omega}](\bar x),\alpha)$ and an arbitrary $\rho>0$.
Set $\rho'=\min\{1,\alpha^{-1}\}\rho$.
By \eqref{strslo} and \eqref{rhoslo}, one can find points $x\in X$ and $\hat\omega=(\omega_1,\ldots,\omega_m)\in \widehat\Omega$ such that $\norm{x-\bar x}<\rho'$, $0<\max_{1\le i\le m} \norm{\omega_i-x}<\rho'$, and
$$
\max\limits_{1\le i\le m} \norm{x-\omega_i}^q-\max\limits_{1\le i\le m} \norm{u-v_i}^q\le \beta\norm{(u,\hat v)-(x,\hat\omega)}_{\rho'}
$$
for all $(u,\hat v)$ with $\hat v=(v_1,\ldots,v_m)\in \widehat\Omega$ near $(x,\hat\omega)$.
In other words, $(x,\hat\omega)$ is a local minimizer of the function
$$
(u,\hat v)\mapsto \max\limits_{1\le i\le m} \norm{u-v_i}^q+\beta\norm{(u,\hat v)-(x,\hat\omega)}_{\rho'}
$$
subject to $\hat v=(v_1,\ldots,v_m)\in \widehat\Omega$.
By definition \eqref{f}, this means that $(x,\hat\omega)$ minimizes locally the function
$$
(u,\hat v)\mapsto f_q(u,\hat v)+\beta\norm{(u,\hat v)-(x,\hat\omega)}_{\rho'},
$$
and consequently its Fr\'echet subdifferential at $(x,\hat\omega)$ contains zero.
Take an
$$
\varepsilon\in\Bigl(0,\min\bigl\{\rho-\norm{x-\bar x},\frac{1}{2}\max_{1\le i\le m} \norm{x-\omega_i},\frac{1}{2}\bigl(\rho-\max_{1\le i\le m} \norm{x-\omega_i}\bigr),\alpha-\beta\bigr\}\Bigr).
$$
Applying Lemma~\ref{l02} and Lemma~\ref{ll01}~(ii), we can find points $x'\in X$, $\hat\omega'=(\omega_1',\ldots,\omega_m')\in \widehat\Omega$, and $(x^*,\hat y^*)\in \partial f_q(x',\hat\omega')$ such that
\begin{gather*}
\norm{x'-x}<\varepsilon,\quad \max_{1\le i\le m} \norm{\omega_i'-\omega_i}<\varepsilon,
\quad\norm{(x^*,\hat y^*)}_{\rho'} =\|x^*\|+\|\hat y^*\|/\rho' <\beta+\varepsilon.
\end{gather*}
It follows that
$$\norm{x'-\bar x}<\rho,\;\; 0<\max_{1\le i\le m} \norm{x'-\omega_i'}<\rho,\;\; \norm{x^*}<\alpha,\;\mbox{ and }\; \norm{\hat y^*}<\rho'\alpha\le \rho.$$
Hence, $\zeta_{\rho,1}^{q*}[\bold{\Omega}](x',\hat\omega')<\alpha$, and consequently $\hat\zeta^{q*}_1[\bold{\Omega}](\bar x)<\alpha$.
By letting $\alpha\to \hat\zeta^q[\bold{\Omega}](\bar x)$, we obtain the claimed inequality.

(ii) follows from (i) and Proposition~\ref{T3.1.1}~(ii).
\qed

\begin{prop}\label{T3.1.3}
Let $X$ be an Asplund space and $\Omega_1$, \ldots, $\Omega_m$ be closed.
\begin{enumerate}
\item
$\hat\zeta^{q*}_2[\bold{\Omega}](\bar x)\le \hat\zeta^{q*}_1[\bold{\Omega}](\bar x)$, where $\hat\zeta^{q*}_1[\bold{\Omega}](\bar x)$ is given by \eqref{strsubslo'},
\begin{equation}\label{c2}
\hat\zeta^{q*}_2[\bold{\Omega}](\bar x):= \lim_{\rho\downarrow 0} \inf_{\substack{\norm{x-\bar x}<\rho\\ \hat\omega=(\omega_1,\ldots,\omega_m)\in \widehat\Omega\\ 0<\max\limits_{1\le i\le m} \norm{x-\omega_i}<\rho}}
\lim_{\varepsilon\downarrow 0} \inf_{\substack{\norm{x'-x}<\varepsilon\\ \hat\omega'\in \widehat\Omega\\ \norm{\hat\omega'-\hat\omega}<\varepsilon}}
\zeta_{\rho,\varepsilon,2}^{q*}[\bold{\Omega}] (x',\hat\omega')
\end{equation}
and, for $x\in X$, $\hat\omega=(\omega_1,\ldots,\omega_m)\in \widehat\Omega$, and $\hat v:=(x-\omega_1,\ldots,x-\omega_m)\ne0$,
\begin{align}
\zeta_{\rho,\varepsilon,2}^{q*}[\bold{\Omega}](x,\hat\omega):= \inf\Biggl\{
\left\|\sum\limits_{i=1}^m{x}_i^*\right\|\mid
&
x_{i}^*\in N_{\Omega_i}(\omega_i) +\rho\B^*\quad (i=1,\ldots,m),
\notag\\&
x_i^*=0
\quad\mbox{if}\quad
\norm{x-\omega_i}<\max_{1\le j\le m} \norm{x-\omega_j},
\notag\\&
\langle x_i^*,x-\omega_i\rangle
\ge\|x_i^*\|(\|x-\omega_i\|-\varepsilon),
\notag\\&
\sum_{i=1}^m\|x_i^*\|={q}\norm{\hat v}^{q-1}
\Biggr\}.\label{c2r}
\end{align}
\item
If $\hat\zeta^{q*}_2[\bold{\Omega}](\bar x)>0$, then $\bold{\Omega}$ is $[q]$-subregular at $\bar x$.
\end{enumerate}
\end{prop}

\textbf{Proof.}
(i)
Let $\rho>0$, $\|x-\bar x\|<\rho$, $\hat\omega:=(\omega_1,\ldots,\omega_m)\in\widehat\Omega$ with $0<\max_{1\le i\le m} \norm{x-\omega_i}<\rho$, $(u^*,\hat v^*) \in\partial f_q(x,\hat\omega)$, where $f_q$ is given by \eqref{f}, and $\|\hat v^*\|<\rho$.
Denote $\hat v:=(x-\omega_1,\ldots,x-\omega_m)$.
Then $0<\|\hat v\|<\rho$.
Observe that function $f_q$ is the sum of two functions on $X^{m+1}$:
$$(x,\hat x)\mapsto \varphi(x,\hat x):= \|(x-x_1,\ldots,x-x_m)\|^q\quad\mbox{and}\quad(x,\hat x) \mapsto\delta_{\widehat\Omega}(\hat x),$$
where $\hat x:=(x_1,\ldots,x_m)$ and $\delta_{\widehat\Omega}$ is the indicator function of $\widehat\Omega$.
The first function is Lipschitz continuous near $(x,\hat\omega)$ (since $\hat v\ne0$), while the second one is lower semicontinuous.
One can apply Lemma~\ref{l02}.
For any $\varepsilon>0$, there exist points $x'\in X$, $\hat x:=(x_1,\ldots,x_m)\in X^m$, $\hat\omega':=(\omega_1',\ldots,\omega_m')\in\widehat\Omega$, $\left(x^*,\hat y^*\right)\in\partial \varphi(x',\hat x)$, and $\hat\omega^*\in N_{\widehat\Omega}(\hat\omega')$ such that
\begin{gather}\notag
\|x'-x\|<\varepsilon,\quad \|\hat x-\hat\omega\|<\frac{\varepsilon}{4},\quad \|\hat\omega'-\hat\omega\|<\frac{\varepsilon}{4},
\\\label{est}
\|(u^*,\hat v^*)-(x^*,\hat y^*)-(0,\hat\omega^*)\|<\varepsilon.
\end{gather}
Taking a smaller $\varepsilon$ if necessary, one can ensure that $\hat v':=(x'-\omega_1',\ldots,x'-\omega_m')\ne0$, $\hat v'':=(x'-x_1,\ldots,x'-x_m)\ne0$, and
\begin{gather}\label{est2}
\|\hat v^*\|+\varepsilon<\rho\left(\frac{\|\hat v'\|}{\|\hat v''\|}\right)^{1-q}
\end{gather}
and, 
for any $i=1,\ldots,m$, $\norm{x'-x_i}<\max_{1\le j\le m} \norm{x'-x_j}$ if and only if $\norm{x'-\omega_i'}<\max_{1\le j\le m} \norm{x'-\omega_j'}$.
By Lemmas~\ref{ll021} and \ref{ll02},
\sloppy
$$\hat x^*:=-\hat y^*
\left(\frac{\|\hat v''\|}{\|\hat v'\|}\right)^{1-q}
\in{q}\norm{\hat v'}^{q-1} J(\hat v'')
\quad\mbox{and}\quad
x^*=x_{1}^*+\ldots+x_{m}^*$$
where
$\hat x^*=(x_1^*,\ldots,x_m^*)$.
By Proposition~\ref{prJ},
\begin{gather*}
\sum_{i=1}^m\|x_i^*\|={q}\norm{\hat v'}^{q-1},
\\
x_i^*=0
\quad\mbox{if}\quad
\norm{x'-\omega_i'}<\max_{1\le j\le m} \norm{x'-\omega_j'},
\end{gather*}
\begin{align*}
\langle x_i^*,x'-\omega_i'\rangle
&\ge\langle x_i^*,x'-x_i\rangle-\|x_i^*\|\,\|x_i-\omega_i'\| =\|x_i^*\|(\|x'-x_i\|-\|x_i-\omega_i'\|)
\\&
\ge\|x_i^*\|(\|x'-\omega_i'\|-2\|x_i-\omega_i'\|)
\ge\|x_i^*\|(\|x'-\omega_i'\|-\varepsilon) \quad(i=1,\ldots,m).
\end{align*}
Inequalities \eqref{est} and \eqref{est2} yield the estimates:
$$
\|u^*\|>\left\|{x}^*\right\|-\varepsilon, \quad \left\|\hat x^*-\hat\omega^*\left(\frac{\|\hat v''\|}{\|\hat v'\|}\right)^{1-q}\right\|<(\|\hat v^*\|+\varepsilon)\left(\frac{\|\hat v''\|}{\|\hat v'\|}\right)^{1-q}<\rho
$$
and consequently
$$
\|u^*\|> \left\|\sum\limits_{i=1}^m{x}_i^*\right\|-\varepsilon, \quad \hat x^*\in N_{\widehat\Omega}(\hat\omega') +\rho\B_m^*.
$$
It follows from Lemma~\ref{ll03} and definitions \eqref{subrhoslo} and \eqref{c2r} that
$$\zeta_{\rho,1}^{q*}[\bold{\Omega}](x,\hat\omega)\ge \zeta_{\rho,\varepsilon,2}^{q*}[\bold{\Omega}] (x',\hat\omega')-\varepsilon.$$
The claimed inequality is a consequence of the last one and definitions \eqref{strsubslo'} and \eqref{c2}.

(ii) follows from (i) and Proposition~\ref{T3.1.2}~(ii).
\qed

\textbf{Proof of Theorem~\ref{T3.1}.}
(i) follows from Proposition~\ref{T3.1.3}~(ii) and definitions \eqref{c2} and \eqref{c2r}.

(ii) follows from \cite[Theorem~4]{Kru09.1} thanks to Remark \ref{compare1}.
\qed


\begin{remark}
One of the main tools in the proof of Theorem~\ref{T3.1} is the fuzzy sum rule (Lemma~\ref{l02}) for Fr\'echet subdifferentials in Asplund spaces.
The statements can be extended to general Banach spaces.
For that, one has to replace Fr\'echet subdifferentials (and normal cones) with some other kind of subdifferentials satisfying a certain set of natural properties including the sum rule (not necessarily fuzzy) -- cf. \cite[p. 345]{KruLop12.1}.

If the sets $\Omega_1$, \ldots $\Omega_m$ are convex or the norm of $X$ is Fr\'echet differentiable away from $0$, then the fuzzy sum rule can be replaced in the proof by either the convex sum rule (Moreau--Rockafellar formula) or the simple (exact) differentiable rule (see, e.g., \cite[Corollary~1.12.2]{Kru03.1}), respectively, to produce dual sufficient conditions for $[q]$-regularity of collections of sets in general Banach spaces in terms of either normals in the sense of convex analysis or Fr\'echet normals.
\end{remark}
\begin{remark}
Since uniform $[q]$-regularity is a stronger property than $[q]$-subregularity (Remark~\ref{rem1}), the criterion in part (ii) of Theorem~\ref{T3.1} is also sufficient for the $[q]$-subregularity (with any $q\in(0,1]$) of the collection of sets in part (i).
\end{remark}

For an example illustrates application of Theorem~\ref{T3.1}~(i) for detecting subregularity of collections of sets, see \cite[Example~4.13]{KruTha13.2}.

\section{$[q]$-regularity of set-valued mappings}\label{S5}

In this section, we present relationships between $[q]$-regularity properties of collections of sets and the corresponding properties of set-valued mappings.
Nonlinear regularity properties of set-valued mappings have been investigated, cf., e.g.,  \cite{GayGeoJea11,YenYaoKie08,Kum09,BorZhuang88,LiMor12,ApeDurStr13,Iof13,FraQui12}.

Consider a set-valued mapping $F:X\rightrightarrows Y$ between metric spaces and a point $(\bar x,\bar y)\in \gr F:=\{(x,y)\in X\times Y\mid y\in F(x)\}$.
\begin{definition}\label{metricR}
\begin{enumerate}
\item
$F$ is \emph{metrically $[q]$-semiregular} at $(\bar{x},\bar{y})$ if there exist positive numbers $\gamma$ and $\delta$ such that
\begin{equation}\label{mhr}
    \gamma d\left(\bar x,F^{-1}(y)\right) \le d^q(y,\bar y),\; \forall y \in B_{\delta}(\bar{y}).
\end{equation}
The exact upper bound of all numbers $\gamma$ such that (\ref{mhr}) is satisfied will be denoted by $\theta^q[F](\bar{x},\bar y)$.
\item
$F$ is \emph{metrically $[q]$-subregular} at $(\bar{x},\bar{y})$ if there exist positive numbers $\gamma$ and $\delta$ such that
\begin{equation}\label{msr}
    \gamma d\left(x,F^{-1}(\bar y)\right) \le d^q(\bar y,F(x)),\; \forall x \in B_{\delta}(\bar{x}).
\end{equation}
The exact upper bound of all numbers $\gamma$ such that (\ref{msr}) is satisfied will be denoted by $\zeta^q[F](\bar{x},\bar y)$.
\item
$F$ is \emph{metrically $[q]$-regular} at $(\bar{x},\bar{y})$ if there exist positive numbers $\gamma$ and $\delta$ such that
\begin{equation}\label{mr}
    \gamma d\left(x,F^{-1}(y)\right) \le d^q\left(y,F(x)\right) ,\; \forall (x,y) \in B_{\delta}(\bar x,\bar{y}).
\end{equation}
The exact upper bound of all numbers $\gamma$ such that (\ref{mr}) is satisfied will be denoted by $\hat{\theta}^q[F](\bar{x},\bar y)$.
\end{enumerate}
\end{definition}

\begin{remark}
Property (ii) and especially property (iii) in Definition \ref{metricR} with $q=1$ are very well known and widely used in variational analysis; see, e.g.,
\cite{RocWet98,Mor06.1,Iof00_,DonRoc09,Kru09.1,Pen89,DonLewRoc03,DonRoc04,ZheNg07,ZheNg10,ZheNg12}. Property (i) (with $q=1$) was introduced in \cite{Kru09.1}.
In \cite{ArtMor11,ApeDurStr13}, it is referred to as \emph{metric hemiregularity}.
\end{remark}

For a collection of sets $\bold{\Omega}:=\{\Omega_1,\ldots,\Omega_m\}$ in a normed linear space $X$, one can consider the set-valued mapping $F:X \rightrightarrows X^m$ defined by (cf. \cite[Proposition~5]{Iof00_}, \cite[Theorem~3]{Kru05.1}, \cite[Proposition 8]{Kru06.1}, \cite[p. 491]{LewLukMal09}, \cite[Proposition 33]{HesLuk})
$$
F(x):= (\Omega_1 -x)\times\ldots\times (\Omega_m -x),\; \forall x \in X.
$$
It is easy to check that, for $x \in X$ and $u=(u_1,\ldots,u_m) \in X^m$, it holds
$$
x \in \bigcap_{i=1}^m\Omega_i \iff 0 \in F(x),\quad F^{-1}(u) = \bigcap_{i=1}^m (\Omega_i-u_i).
$$

The next proposition is a consequence of Theorem~\ref{T2.1}.

\begin{proposition}\label{theorem13}
Consider $\bold{\Omega}$ and $F$ as above and a point $\bar x\in \bigcap_{i=1}^m\Omega_i$.
\begin{enumerate}
\item
$\bold{\Omega}$ is $[q]$-semiregular at $\bar x$ if and only if $F$ is metrically $[q]$-semiregular at $(\bar x,0)$.
Moreover, $\theta^q[\bold{\Omega}](\bar{x})=\theta^q[F](\bar{x},0)$.
\item
$\bold{\Omega}$ is $[q]$-subregular at $\bar x$ if and only if $F$ is metrically $[q]$-subregular at $(\bar x,0)$.
Moreover, $\zeta^q[\bold{\Omega}](\bar{x})=\zeta^q[F](\bar{x},0)$.
\item
$\bold{\Omega}$ is uniformly $[q]$-regular at $\bar x$ if and only if $F$ is metrically $[q]$-regular at $(\bar x,0)$.
Moreover, $\hat{\theta}^q[\bold{\Omega}](\bar{x})=\hat{\theta}^q[F](\bar{x},0)$.
\end{enumerate}
\end{proposition}

For a further discussion of the relationships between regularity properties of $\bold{\Omega}$ and $F$ see \cite[Remark 5.4]{KruTha13.2}.

Conversely, regularity properties of set-valued mappings between normed linear spaces can be treated as realizations of the corresponding properties of certain collections of two sets.

For a given set-valued mapping $F:X\rightrightarrows Y$ between normed linear spaces and a point $(\bar x,\bar y)\in \gr F$, one can consider the collection $\bold{\Omega}$ of two sets $\Omega_1= \gr F$ and $\Omega_2=X\times \{\bar y\}$ in $X\times Y$.
It is clear that $(\bar x,\bar y)\in \Omega_1\cap\Omega_2$.
\begin{proposition}\label{theorem14}
Consider $F$ and $\bold{\Omega}$ as above.
\begin{enumerate}
\item
$F$ is metrically $[q]$-semiregular at $(\bar x,\bar y)$ if and only if $\bold{\Omega}$ is $[q]$-semiregular at $(\bar x,\bar y)$. Moreover,
\begin{equation}\label{newestSem}
\frac{\theta^q[F](\bar{x},\bar y)}{\theta^q[F](\bar{x},\bar y)+2^q}\le \theta^q[\bold{\Omega}](\bar{x},\bar y)\le \theta^q[F](\bar{x},\bar y)/2^q.
\end{equation}
\item
$F$ is metrically $[q]$-subregular at $(\bar x,\bar y)$ if and only if $\bold{\Omega}$ is $[q]$-subregular at $(\bar x,\bar y)$. Moreover,
\begin{equation}\label{newest}
\frac{\zeta^q[F](\bar{x},\bar y)}{\zeta^q[F](\bar{x},\bar y)+2^q}\le \zeta^q[\bold{\Omega}](\bar{x},\bar y)\le \zeta^q[F](\bar{x},\bar y)/2^q.
\end{equation}
\item
$F$ is metrically $[q]$-regular at $(\bar x,\bar y)$ if and only if $\bold{\Omega}$ is uniformly $[q]$-regular at $(\bar x,\bar y)$. Moreover,
\begin{equation}\label{newestReg}
\frac{\hat{\theta}^q[F](\bar{x},\bar y)}{\hat{\theta}^q[F](\bar{x},\bar y)+2^q}\le \hat{\theta}^q[\bold{\Omega}](\bar{x},\bar y)\le \hat{\theta}^q[F](\bar{x},\bar y)/2^q.
\end{equation}
\end{enumerate}
\end{proposition}
\textbf{Proof.} (i) Suppose $F$ is metrically $[q]$-semiregular at $(\bar x,\bar y)$, i.e., $\theta^q[F](\bar{x},\bar y)>0$. Fix a $\gamma\in(0,\theta^q[F](\bar{x},\bar y))$. Then there exists a number $\delta'>0$ such that \eqref{mhr} is satisfied for all $y\in B_{\delta'}(\bar y)$. Set an $\alpha:=\frac{\gamma}{\gamma+2^q}$ (so $2^q\alpha/\gamma+\alpha^{\frac{1}{q}}<1$) and a $\delta:=\min\left\{\frac{\delta'^q}{2^q\alpha},1\right\}$. We are going to check that
\begin{equation}\label{HR'}
\left(\Omega_1-(u_1,v_1)\right)\bigcap \left(\Omega_2-(u_2,v_2)\right)\bigcap B_{\rho}(\bar x,\bar y)\neq \emptyset
\end{equation}
for all $\rho\in (0,\delta)$ and $(u_1,v_1),(u_2,v_2)\in (\alpha\rho)^{\frac{1}{q}}\mathbb{B}$.
Indeed, take any $\rho\in (0,\delta)$ and $(u_1,v_1),(u_2,v_2)\in (\alpha\rho)^{\frac{1}{q}}\mathbb{B}$. We need to find a point $(x,y)\in B_{\rho}(\bar x,\bar y)$ satisfying
$$
\left\{
  \begin{array}{ll}
    (x,y)+(u_1,v_1)\in \gr F,\\
    y=\bar y-v_2\hbox{.}
  \end{array}
\right.
$$
We set $y':=\bar y-v_2+v_1$, so $y' \in B_{\delta'}(\bar y)$ as $\|y'-\bar y\|=\|v_1-v_2\|\le 2(\alpha\rho)^{\frac{1}{q}}<2(\alpha\delta)^{\frac{1}{q}}=\delta'$. Then, by \eqref{mhr}, there is an $x'\in F^{-1}(y')$ such that
$$
\|\bar x-x'\|\le \frac{1}{\gamma}\|\bar y-y'\|^q.
$$
Put $y:=y'-v_1=\bar y-v_2$ and $x:=x'-u_1$. Then it holds
$(x,y)+(u_1,v_1)=(x',y') \in \gr F$, $\norm{y-\bar y}=\|v_2\|\le (\alpha\rho)^{\frac{1}{q}}<\rho$, and
\begin{multline*}
\|x-\bar x\|\le \|x-x'\|+\|x'-\bar x\|\le \|u_1\|+\frac{1}{\gamma}\|\bar y-y'\|^q\\
=\|u_1\|+\frac{1}{\gamma}\|v_1-v_2\|^q \le (2^q\alpha/\gamma+\alpha^{\frac{1}{q}})\rho<\rho.
\end{multline*}
Hence, \eqref{HR'} is proved.

The above reasoning also yields the first inequality in \eqref{newestSem}.

To prove the inverse implication, we suppose $\bold{\Omega}$ is $[q]$-semiregular at $(\bar x,\bar y)$, i.e., $\theta^q[\bold{\Omega}](\bar{x},\bar y)>0$. Fix an $\alpha\in (0,\theta^q[\bold{\Omega}](\bar{x},\bar y))$. Then there exists $\delta'>0$ such that \eqref{HR'} holds true for all $\rho\in (0,\delta')$ and $(u_1,v_1),(u_2,v_2)\in (\alpha\rho)^{\frac{1}{q}}\mathbb{B}$.
Set $\gamma:=2^q\alpha$ and $\delta<(\alpha\delta')^{\frac{1}{q}}$. We are going to check that \eqref{mhr} is satisfied.
Take any $y\in B_{\delta}(\bar y)$, i.e., $\|y-\bar y\|\le \delta<(\alpha\delta')^{\frac{1}{q}}$. Set $r\in (0,\delta')$ such that $\|y-\bar y\|=(\alpha r)^{\frac{1}{q}}$.
Then, applying \eqref{HR'} for $\rho:=\frac{r}{2^q}\in (0,\delta')$, and $(u_1,v_1):=\left(0,\frac{y-\bar y}{2}\right)$, $(u_2,v_2):=\left(0,\frac{\bar y-y}{2}\right)\in \left(\alpha\frac{r}{2^q}\right)^{\frac{1}{q}}\mathbb{B}$, we can find
$(x_1,y_1)\in \gr F$ and $(x_2,\bar y)\in \Omega_2$ satisfying
$$
(x_1,y_1)-(u_1,v_1)=(x_2,\bar y)-(u_2,v_2)\in B_{\frac{r}{2^q}}(\bar x,\bar y).
$$
This implies that $y_1=y$, $x_1\in F^{-1}(y)$, and
$$
\|x_1-\bar x\|\le \frac{r}{2^q}=\frac{1}{2^q\alpha}\|y-\bar y\|^q=\frac{1}{\gamma}\|y-\bar y\|^q.
$$
Hence, \eqref{mhr} holds true.

The last reasoning also yields the second inequality in \eqref{newestSem}.

(ii) Suppose $F$ is metrically $[q]$-subregular at $(\bar x,\bar y)$, i.e., $\zeta^q[F](\bar{x},\bar y)>0$. Fix a $\gamma\in (0,\zeta^q[F](\bar{x},\bar y))$. Then there exists a $\delta'>0$ (one can take $\delta'\in (0,1)$) such that \eqref{msr} is satisfied for all $x\in B_{\delta'}(\bar x)$.
Set an $\alpha:=\frac{\gamma}{\gamma+2^q}$ (so $2^q\alpha/\gamma+\alpha^{\frac{1}{q}}<1$) and a $\delta>0$ satisfying $(\alpha\delta)^{\frac{1}{q}}+\delta<\delta'$. We are going to check that
\begin{equation}\label{Hlr''}
\left(\Omega_1+ (\alpha\rho)^{\frac{1}{q}}\mathbb{B}\right)\bigcap \left(\Omega_2+ (\alpha\rho)^{\frac{1}{q}}\mathbb{B}\right)\bigcap B_{\delta}(\bar x,\bar y)\subseteq \Omega_1\cap \Omega_2+\rho\mathbb{B}
\end{equation}
for all $\rho\in (0,\delta)$.
Indeed, take any
$$
(x,y)\in\left(\Omega_1+ (\alpha\rho)^{\frac{1}{q}}\mathbb{B}\right)\bigcap \left(\Omega_2+ (\alpha\rho)^{\frac{1}{q}}\mathbb{B}\right)\bigcap B_{\delta}(\bar x,\bar y).
$$
Then
$(x,y)=(x_1,y_1)+(u_1,v_1)=(x_2,\bar y)+(u_2,v_2)$ for some $(x_1,y_1)\in \gr F$, $x_2\in X$, and $(u_1,v_1)$, $(u_2,v_2)\in (\alpha\rho)^{\frac{1}{q}}\mathbb{B}$.
Since
$$\|x_1-\bar x\|\le\|u_1\|+\|x-\bar x\| \le (\alpha\rho)^{\frac{1}{q}}+\delta<\delta',$$
by \eqref{msr}, there exists an $x'\in F^{-1}(\bar y)$ such that $\|x_1-x'\|\le \frac{1}{\gamma}\|\bar y-y_1\|^q$. Then
\begin{align*}
\norm{x_1-x'+u_1}\le
&\frac{1}{\gamma}\|\bar y-y_1\|^q+\norm{u_1}=\frac{1}{\gamma}\norm{v_1-v_2}^q+\norm{u_1}
\\
\le &\frac{2^q\alpha\rho}{\gamma}+(\alpha\rho)^{\frac{1}{q}}\le \Big{(}\frac{2^q\alpha}{\gamma}+\alpha^{\frac{1}{q}}\Big{)}\rho<\rho,
\\
\norm{v_2}\le &(\alpha\rho)^{\frac{1}{q}}\le \alpha^{\frac{1}{q}}\rho<\rho.
\end{align*}
Hence, $(x,y)=(x',\bar y)+(x_1-x'+u_1,v_2)\in \Omega_1\cap\Omega_2+\rho\mathbb{B}$.

The above reasoning also yields the first inequality in \eqref{newest}.

To prove the inverse implication, we suppose that $\bold{\Omega}$ is $[q]$-subregular at $(\bar x,\bar y)$, i.e., $\zeta^q[\bold{\Omega}](\bar{x},\bar y)>0$. Fix an $\alpha\in (0,\zeta^q[\bold{\Omega}](\bar{x},\bar y))$. Then there exists a $\delta'>0$ such that \eqref{Hlr''} holds true for all $\rho\in (0,\delta')$.
Set $\gamma:=2^q\alpha>0$ and $\delta:=\min\left\{\delta',\gamma\delta',\frac{2^q\delta'^q}{\gamma}\right\}$. We are going to check that \eqref{msr} holds true.
Take any $x\in B_{\delta}(\bar x)$.
Because $d(x,F^{-1}(\bar y))\le \|x-\bar x\|\le \delta$, it is sufficient to consider the case $0<d(\bar y,F(x))<(\gamma\delta)^{\frac{1}{q}}$. We take a $y\in F(x)$ such that $d(\bar y, F(x))\le \|y-\bar y\|:=r <(\gamma\delta)^{\frac{1}{q}}$. Then
$$
\left(x,\frac{y+\bar y}{2}\right)=(x,y)+\left(0,\frac{\bar y-y}{2}\right)=(x,\bar y)+\left(0,\frac{y-\bar y}{2}\right),\quad\norm{\frac{\bar y-y}{2}}= \frac{ r}{2}<\delta',
$$
and consequently
\begin{equation}\label{e43}
\left(x,\frac{y+\bar y}{2}\right)\in \left(\Omega_1+ \frac{ r}{2}\mathbb{B}\right)\bigcap \left(\Omega_2+ \frac{ r}{2}\mathbb{B}\right)\bigcap B_{\delta'}(\bar x,\bar y).
\end{equation}
Take $\rho:=\frac{r^q}{2^q\alpha}< \delta\le \delta'$.
Then $\frac{ r}{2}=(\alpha\rho)^{\frac{1}{q}}$, and it follows from \eqref{Hlr''} and \eqref{e43} that
$$
\left(x,\frac{y+\bar y}{2}\right)\in \Omega_1\cap \Omega_2+\frac{ r^q}{2^q\alpha}\mathbb{B} =F^{-1}(\bar y)\times\{\bar y\} +\frac{\|y-\bar y\|^q}{\gamma}\mathbb{B}.
$$
Hence, there is an $x'\in F^{-1}(\bar y)$ such that
$$
\|x-x'\|\le \frac{1}{\gamma}\|y-\bar y\|^q.
$$
Taking infimum in the last inequality over $x'\in F^{-1}(\bar y)$ and $y\in F(x)$, we arrive at \eqref{msr}.

(iii) Suppose $F$ is metrically $[q]$-regular at $(\bar x,\bar y)$, i.e., $\hat{\theta}^q[F](\bar{x},\bar y)>0$. Fix a $\gamma\in (0,\hat{\theta}^q[F](\bar{x},\bar y))$. Then there exists a $\delta'>0$ (one can take $\delta'\in (0,1)$) such that \eqref{mr} is satisfied for all $(x,y)\in B_{\delta'}(\bar x,\bar y)$.
Set an $\alpha:=\frac{\gamma}{\gamma+2^q}$ (so $2^q\alpha/\gamma+\alpha^{\frac{1}{q}}<1$) and a $\delta:=\frac{\delta'}{2\alpha^{\frac{1}{q}}+1}$. We are going to check that
\begin{equation}\label{HUR'}
\left(\Omega_1-(x_1,y_1)-(u_1,v_1)\right)\bigcap \left(\Omega_2-(x_2,\bar y)-(u_2,v_2)\right)\bigcap (\rho\mathbb{B})\neq \emptyset
\end{equation}
for all $\rho\in (0,\delta)$, $(x_1,y_1)\in \Omega_1\cap B_{\delta}(\bar x,\bar y),x_2\in B_{\delta}(\bar x)$, and $(u_1,v_1),(u_2,v_2)\in (\alpha\rho)^{\frac{1}{q}}\mathbb{B}$.
Take any such $\rho,(x_1,y_1),x_2,(u_1,v_1)$, and $(u_2,v_2)$. We need to find $(a,b)\in \rho\mathbb{B}$ satisfying
$$
\left\{
  \begin{array}{ll}
    (x_1,y_1)+(u_1,v_1)+(a,b)\in \gr F,\\
    b=-v_2\hbox{.}
  \end{array}
\right.
$$
We set $y'=y_1-v_2+v_1$, so $y' \in B_{\delta'}(\bar y)$ as
$$
\|y'-\bar y\|\leq \|y'-y_1\|+\|y_1-\bar y\|\le \|v_1-v_2\|+\delta \le 2(\alpha\rho)^{\frac{1}{q}}+\delta<(2\alpha^{\frac{1}{q}}+1)\delta=\delta'.
$$
Then, applying \eqref{mr} for $(x_1,y')\in B_{\delta'}(\bar x,\bar y)$, we find $x'\in F^{-1}(y')$ such that
$$
\|x_1-x'\|\le \frac{1}{\gamma}d^q(y',F(x_1))\le\frac{1}{\gamma}\|y'-y_1\|^q=\frac{1}{\gamma}\|v_1-v_2\|^q\le \frac{2^q\alpha\rho}{\gamma} .
$$
Put $a=x'-x_1-u_1$ and $b=-v_2$. Then $\|a\|\le \|x'-x_1\|+\|u_1\|\le (2^q\alpha/\gamma+\alpha^{\frac{1}{q}})\rho<\rho$, $\|b\|\le (\alpha\rho)^{\frac{1}{q}}<\rho$, and it holds
$(x_1,y_1)+(u_1,v_1)+(a,b)=(x',y') \in \gr F$.\\
Hence, \eqref{HUR'} is proved.

The above reasoning also yields the first inequality in \eqref{newestReg}.

To prove the inverse implication, we suppose that $\bold{\Omega}$ is uniformly $[q]$-regular at $(\bar x,\bar y)$, i.e., $\hat{\theta}^q[\bold{\Omega}](\bar{x},\bar y)>0$. Fix an $\alpha\in (0,\hat{\theta}^q[\bold{\Omega}](\bar{x},\bar y))$. Then there exists a $\delta'>0$ (one can take $\delta'\in (0,1)$) such that \eqref{HUR'} holds true for all $\rho\in (0,\delta')$, $(x_1,y_1)\in \Omega_1\cap B_{\delta'}(\bar x,\bar y),x_2\in B_{\delta'}(\bar x)$, and $(u_1,v_1),(u_2,v_2)\in (\alpha\rho)^{\frac{1}{q}}\mathbb{B}$.
Set $\gamma:=2^q\alpha>0$. Because $\theta^q[\bold{\Omega}](\bar{x},\bar y)\ge \hat{\theta}^q[\bold{\Omega}](\bar{x},\bar y)$ (see Remark \ref{rem1}), assertion (i) implies that there exists a $\delta^*>0$ such that \eqref{mhr} is satisfied for all $y\in B_{\delta^*}(\bar y)$.
Choose a positive number $\delta$ satisfying the following conditions
\begin{equation}\label{del}
\left\{
  \begin{array}{ll}
    \delta\le \delta^*,\\
    2^q\delta+\frac{\delta^q}{\alpha}\le \delta',\\
    (2^q\alpha\delta+\delta^q)^{\frac{1}{q}}+\delta \le \delta'\hbox{.}
  \end{array}
\right.
\end{equation}
Now, take any $(x,y)\in B_{\delta}(\bar x,\bar y)$.
We are going to check that \eqref{mr} is satisfied. Because \eqref{mhr} implies
\begin{equation*}
    \gamma d(x,F^{-1}(y))\le \gamma\|x-\bar x\|+\gamma d(\bar x,F^{-1}(y))\le \gamma\delta+\|y-\bar y\|^q\le \gamma\delta+\delta^q,
\end{equation*}
it suffices to consider the case $d(y,F(x))<(\gamma\delta+\delta^q)^{\frac{1}{q}}$ (note that $\gamma\delta+\delta^q\le \alpha\delta'$ by \eqref{del}.)
Choose a $y'\in F(x)$ such that
$$d(y,F(x))\le \|y-y'\|<(\gamma\delta+\delta^q)^{\frac{1}{q}}$$
and set $r\in (0,\delta')$ such that $\|y-y'\|=(\alpha r)^{\frac{1}{q}}$.
Then
$$\|y'-\bar y\|\le \|y'-y\|+\|y-\bar y\|<(\gamma\delta+\delta^q)^{\frac{1}{q}}+\delta\le \delta'$$
due to \eqref{del}.
Applying \eqref{HUR'} with
$$(x_1,y_1):=(x,y')\in \gr F\cap B_{\delta'}(\bar x,\bar y),\quad (x_2,y_2):=(\bar x,\bar y),$$ $$(u_1,v_1):=\left(0,\frac{y-y'}{2}\right),\quad (u_2,v_2):=\left(0,\frac{y'-y}{2}\right)\in \left(\alpha\frac{r}{2^q}\right)^{\frac{1}{q}}\mathbb{B},$$ we can find
$(\tilde x,\tilde y)\in \gr F$ and $(z,\bar y)\in \Omega_2$ satisfying
$$
(\tilde x,\tilde y)-(x_1,y_1)-(u_1,v_1)=(z,\bar y)-(x_2,\bar y)-(u_2,v_2)\in \frac{r}{2^q}\mathbb{B}.
$$
This implies $\tilde x-x_1\in \frac{r}{2^q}\mathbb{B}$ and $\tilde y=y_1+v_1-v_2=y$, so $\tilde x\in F^{-1}(y)$. Then we obtain
$$
d(x,F^{-1}(y))\le \|x-\tilde x\|\le \frac{r}{2^q}=\frac{1}{2^q\alpha}\|y-y'\|^q=\frac{1}{\gamma}\|y-y'\|^q.
$$
Taking infimum in the last inequality over $y'\in F(x)$, we arrive at \eqref{mr}.

The last reasoning also yields the second inequality in \eqref{newestReg}.
\qed

{\bf Acknowledgements}.
The authors wish to thank the anonymous referees for very careful reading of the paper and many valuable comments and suggestions, which helped us improve the presentation.

The research was supported by the Australian Research Council, project DP110102011.

\begin{thebibliography}{10}
\expandafter\ifx\csname url\endcsname\relax
  \def\url#1{\texttt{#1}}\fi
\expandafter\ifx\csname
urlprefix\endcsname\relax\def\urlprefix{URL }\fi

\bibitem{AnhKruTha}
L.~Q. Anh, A.~Y. Kruger, N.~H. Thao, On H\"older calmness of
solution mappings in parametric equilibrium problems, TOP, DOI:
10.1007/s11750-012-0259-3.

\bibitem{ApeDurStr13}
M. Apetrii, M. Durea, R. Strugariu, On subregularity properties
of set-valued mappings, Set-Valued Var. Anal. 21~(1) (2013) 93--126.

\bibitem{ArtMor11}
F.~J. Arag\'on Artacho, B.~S. Mordukhovich, Enhanced metric
regularity and Lipschitzian properties of variational systems, J. Global
Optim. 50~(1) (2011) 145--167.

\bibitem{AttBolRedSou10}
H. Attouch, J. Bolte, P. Redont, A. Soubeyran, Proximal
alternating minimization and projection methods for nonconvex problems:
an approach based on the Kurdyka-Lojasiewicz inequality, Math. Oper.
Res. 35~(2) (2010) 438--457.

\bibitem{AusDanThi05}
D. Aussel, A. Daniilidis, L. Thibault, Subsmooth sets: functional
characterizations and related concepts, Trans. Amer. Math. Soc.
357~(4) (2005) 1275--1301.

\bibitem{BakDeuLi05}
A. Bakan, F. Deutsch, W. Li, Strong CHIP, normality, and
linear regularity of convex sets, Trans. Amer. Math. Soc. 357~(10) (2005) 3831--3863.

\bibitem{BauBor93}
H.~H. Bauschke, J.~M. Borwein, On the convergence of von
Neumann's alternating projection algorithm for two sets, Set-Valued Anal.
1~(2) (1993) 185--212.

\bibitem{BauBor96}
H.~H. Bauschke, J.~M. Borwein, On projection algorithms for
solving convex feasibility problems, SIAM Rev. 38~(3) (1996) 367--426.

\bibitem{BauBorLi99}
H.~H. Bauschke, J.~M. Borwein, W. Li, Strong conical hull
intersection property, bounded linear regularity, Jameson's property (G),
and error bounds in convex optimization, Math. Program., Ser. A 86~(1) (1999) 135--160.

\bibitem{BauBorTse00}
H.~H. Bauschke, J.~M. Borwein, P. Tseng, Bounded linear
regularity, strong CHIP, and CHIP are distinct properties, J. Convex Anal.
7~(2) (2000) 395--412.

\bibitem{BorZhuang88}
J.~M. Borwein, D.~M. Zhuang, Verifiable necessary and suffcient conditions for openness and regularity for set-valued and single-valued
mapps, J. Math. Anal. Appl. 134 (1988) 441--459.

\bibitem{BurDeng05}
J.~V. Burke, S. Deng, Weak sharp minima revisited. II. Application
to linear regularity and error bounds, Math. Program., Ser. B 104~(2-3) (2005)
235--261.

\bibitem{DonLewRoc03}
A.~L. Dontchev, A.~S. Lewis, R.~T. Rockafellar, The radius
of metric regularity, Trans. Amer. Math. Soc. 355~(2) (2003) 493--517.

\bibitem{DonRoc04}
A.~L. Dontchev, R.~T. Rockafellar, Regularity and conditioning
of solution mappings in variational analysis, Set-Valued Anal. 12~(1-2) (2004) 79--109.

\bibitem{DonRoc09}
A.~L. Dontchev, R.~T. Rockafellar, Implicit Functions and Solution Mappings. A View from Variational Analysis. Springer Monographs
in Mathematics, Springer, Dordrecht, 2009.

\bibitem{Eke74}
I. Ekeland, On the variational principle. J. Math. Anal. Appl. 47 (1974) 324--353.

\bibitem{Fab89}
M. Fabian, Subdifferentiability and trustworthiness in the light of a new variational principle of {B}orwein and {P}reiss. Acta Univ. Carolinae 30 (1989) 51--56.

\bibitem{Fra89}
H. Frankowska, High order inverse function theorems, Ann. Inst. H.
Poincar\'e Anal. Non Lin\'eaire 6, suppl. (1989) 283--303.

\bibitem{FraQui12}
H. Frankowska, M. Quincampoix, H\"older metric regularity of
set-valued maps, Math. Program., Ser. A 132~(1-2) (2012) 333--354.

\bibitem{GayGeoJea11}
M. Gaydu, M.~H. Geoffroy, C. Jean-Alexis, Metric subregularity
of order $q$ and the solving of inclusions, Cent. Eur. J. Math. 9~(1)
(2011) 147--161.

\bibitem{HesLuk}
R. Hesse, D. R. Luke, Nonconvex notions of regularity and convergence of fundamental algorithms for feasibility problems. Preprint, arXiv:1212.3349v1.
Accessed 17/12/2013.

\bibitem{Hua12.1}
X.~X. Huang, Calmness and exact penalization in constrained scalar
set-valued optimization, J. Optim. Theory Appl. 154~(1) (2012) 108--119.

\bibitem{Iof89}
A.~D. Ioffe, Approximate subdifferentials and applications. III. The metric
theory, Mathematika 36~(1) (1989) 1--38.

\bibitem{Iof00_}
A.~D. Ioffe, Metric regularity and subdifferential calculus, Russian Math. Surveys 55 (2000) 501--558.

\bibitem{Iof13}
A.~D. Ioffe, Nonlinear regularity models, Math. Program. 139~(1-2) (2013) 223--242.

\bibitem{KlaLi99}
D. Klatte, W. Li, Asymptotic constraint qualifications and global
error bounds for convex inequalities, Math. Program., Ser. A 84~(1) (1999) 137--160.

\bibitem{Kru96.2}
A.~Y. Kruger, Strict $\varepsilon$-semidifferentials and differentiation of multivalued
mappings, Dokl. Akad. Nauk Belarusi 40~(6) (1996) 38--43, in Russian.

\bibitem{Kru98.1}
A.~Y. Kruger, On the extremality of set systems, Dokl. Nats. Akad.
Nauk Belarusi 42~(1) (1998) 24--28, in Russian.

\bibitem{Kru00.1}
A.~Y. Kruger, Strict $(\varepsilon,\delta)$-semidifferentials and the extremality of sets
and functions, Dokl. Nats. Akad. Nauk Belarusi 44~(2) (2000) 19--22, in Russian.

\bibitem{Kru02.1}
A.~Y. Kruger, Strict $(\varepsilon,\delta)$-subdifferentials and extremality conditions,
Optimization 51~(3) (2002) 539--554.

\bibitem{Kru03.1}
A.~Y. Kruger, On Fr\'echet subdifferentials, J. Math. Sci. 116~(3) (2003) 3325--3358.

\bibitem{Kru04.1}
A.~Y. Kruger, Weak stationarity: eliminating the gap between necessary
and sufficient conditions, Optimization 53~(2) (2004) 147--164.

\bibitem{Kru05.1}
A.~Y. Kruger, Stationarity and regularity of set systems, Pac. J. Optim.
1~(1) (2005) 101--126.

\bibitem{Kru06.1}
A.~Y. Kruger, About regularity of collections of sets, Set-Valued Anal.
14~(2) (2006) 187--206.

\bibitem{Kru09.1}
A.~Y. Kruger, About stationarity and regularity in variational analysis,
Taiwanese J. Math. 13(6A) (2009) 1737--1785.

\bibitem{KruLop12.1}
A.~Y. Kruger, M.~A. L\'opez, Stationarity and regularity of infinite
collections of sets, J. Optim. Theory Appl. 154~(2) (2012) 339--369.

\bibitem{KruLop12.2}
A.~Y. Kruger, M.~A. L\'opez, Stationarity and regularity of infinite
collections of sets. Applications to infinitely constrained optimization, J.
Optim. Theory Appl. 155~(2) (2012) 390--416.

\bibitem{KruTha13.1}
A.~Y. Kruger, N.~H. Thao, About uniform regularity of collections
of sets, Serdica Math. J. 39 (2013) 287--312.

\bibitem{KruTha13.2}
A.~Y. Kruger, N.~H. Thao, Quantitative characterizations of regularity properties of collections of sets, to be published.

\bibitem{Kum09}
B. Kummer, Inclusions in general spaces: Hoelder stability, solution
schemes and Ekeland's principle, J. Math. Anal. Appl. 358~(2) (2009) 327--344.

\bibitem{LewLukMal09}
A.~S. Lewis, D.~R. Luke, J. Malick, Local linear convergence for
alternating and averaged nonconvex projections, Found. Comput. Math. 9~(4) (2009)
485--513.

\bibitem{LewPan98}
A.~S. Lewis, J.-S. Pang, Error bounds for convex inequality systems.
In Generalized Convexity, Generalized Monotonicity: Recent Results
(Luminy, 1996), Kluwer Acad. Publ., Dordrecht (1998) 75--110.

\bibitem{LiNgPon07}
C. Li, K.~F. Ng, T.~K. Pong, The SECQ, linear regularity, and
the strong CHIP for an infinite system of closed convex sets in normed
linear spaces, SIAM J. Optim. 18~(2) (2007) 643--665.

\bibitem{LiMor12}
G. Li, B.~S. Mordukhovich, H\"older metric subregularity with
applications to proximal point method, SIAM J. Optim. 22~(4) (2012) 1655--1684.

\bibitem{LiNahSin00}
W. Li, C. Nahak, I. Singer, Constraint qualifications for semiinfinite systems
 of convex inequalities, SIAM J. Optim. 11~(1) (2000) 31--52.

\bibitem{Lucc06}
R. Lucchetti, Convexity and Well-Posed Problems. CMS Books in Mathematics/
Ouvrages de Math\'ematiques de la SMC, Springer, New York, 2006.

\bibitem{Luk12}
D.~R. Luke, Local linear convergence of approximate projections onto
regularized sets, Nonlinear Anal. 75~(3) (2012) 1531--1546.

\bibitem{Luk13}
D.~R. Luke, Prox-regularity of rank constraint sets and implications for algorithms, J. Math. Imaging Vis. (2013), DOI 10.1007/s10851-012-0406-3.

\bibitem{Mor06.1}
B.~S. Mordukhovich, Variational Analysis and Generalized Differentiation. I: Basic Theory, Springer-Verlag, Berlin, 2006.

\bibitem{NgaiThe01}
H.~V. Ngai, M. Th\'era, Metric inequality, subdifferential calculus
and applications, Set-Valued Anal. 9~(1-2) (2001) 187--216.

\bibitem{Pen89}
J.-P. Penot, Metric regularity, openness and Lipschitz behavior of multifunctions,
Nonlinear Anal. 13 (1989) 629--643.

\bibitem{Pen13}
J.-P. Penot, Calculus Without Derivatives, Springer-Verlag, New York, 2013.

\bibitem{Phelps}
R.~R. Phelps, \emph{{Convex Functions, Monotone Operators and
  Differentiability, 2nd edition}}, Lecture Notes in Mathematics, Vol. 1364,
  Springer, New York, 1993.

\bibitem{RocWet98}
R.~T. Rockafellar, R.~J.-B. Wets, Variational Analysis,
Springer-Verlag, Berlin, 1998.

\bibitem{YenYaoKie08}
N.~D. Yen, J.-C. Yao, B.~T. Kien, Covering properties at
positive-order rates of multifunctions and some related topics, J. Math.
Anal. Appl. 338~(1) (2008) 467--478.

\bibitem{Zal02}
C. Z{\u{a}}linescu, Convex Analysis in General Vector Spaces. World Scientific Publishing Co. Inc.,
River Edge, NJ, 2002.

\bibitem{ZheNg07}
X.~Y. Zheng, K.~F. Ng, Metric subregularity and constraint qualifications
for convex generalized equations in Banach spaces, SIAM J. Optim.
18 (2007) 437--460.

\bibitem{ZheNg08}
X.~Y. Zheng, K.~F. Ng, Linear regularity for a collection of subsmooth
sets in Banach spaces, SIAM J. Optim. 19~(1) (2008) 62--76.

\bibitem{ZheNg10}
X.~Y. Zheng, K.~F. Ng, Metric subregularity and calmness for
nonconvex generalized equations in Banach spaces, SIAM J. Optim. 20~(5) (2010) 2119--2136.

\bibitem{ZheNg12}
X.~Y. Zheng, K.~F. Ng, Metric subregularity for proximal generalized
equations in Hilbert spaces, Nonlinear Anal. 75~(3) (2012) 1686--1699.

\bibitem{ZheWeiYao10}
X.~Y. Zheng, Z. Wei, J.-C. Yao, Uniform subsmoothness and
linear regularity for a collection of infinitely many closed sets, Nonlinear
Anal. 73~(2) (2010) 413--430.

\end{thebibliography}

\end{document}